\documentclass[pdflatex,sn-mathphys-num]{sn-jnl}

\usepackage{graphicx}
\usepackage{amsmath,amssymb}
\usepackage{amsthm}
\usepackage{xcolor}
\usepackage{booktabs}
\usepackage{bm}
\usepackage{epstopdf}

\usepackage{tikz}
\usetikzlibrary{shapes, patterns}

\usepackage{siunitx}

\usepackage[capitalise]{cleveref}
\usepackage{doi}

\newcommand{\pvct}[1]{\bm{#1}}
\newcommand{\vct}[1]{\bm{\mathsf{#1}}}
\newcommand{\mtx}[1]{\bm{\mathsf{#1}}}

\newcommand{\pxx}{\pvct{x}}
\newcommand{\hpxx}{\hat{\pvct{x}}}
\newcommand{\pyy}{\pvct{y}}
\newcommand{\pzz}{\pvct{z}}

\definecolor{UTburnt}{RGB}{191,87,0}
\definecolor{UTorange}{RGB}{248,151,31}
\definecolor{UTyellow}{RGB}{255,214,0}
\definecolor{UTpea}{RGB}{166,205,87}
\definecolor{UTgreen}{RGB}{87,157,66}
\definecolor{UTcyan}{RGB}{0,169,183}
\definecolor{UTblue}{RGB}{0,95,134}
\definecolor{UTmetal}{RGB}{156,173,183}
\definecolor{UTbeige}{RGB}{214,210,196}
\definecolor{UTdark}{RGB}{51,63,72}

\definecolor{myred}{RGB}{231,51,57}
\definecolor{mycrimson}{RGB}{176,9,48}
\definecolor{mypink}{RGB}{219,81,151}
\colorlet{mypurple}{UTblue!50!mycrimson}

\theoremstyle{thmstyleone}

\theoremstyle{thmstyletwo}
\newtheorem{example}{Example}

\theoremstyle{thmstylethree}

\begin{document}

\title{A stable and fast method for solving multibody scattering
problems via the method of fundamental solutions}

\author*[1]{\fnm{Yunhui} \sur{Cai}}\email{caiyh@utexas.edu}

\author[2]{\fnm{Joar} \sur{Bagge}}\email{joar.bagge@austin.utexas.edu}

\author[1,2]{\fnm{Per-Gunnar} \sur{Martinsson}}\email{pgm@oden.utexas.edu}

\affil[1]{\orgdiv{Department of Mathematics}, \orgname{The University of Texas at Austin}, \orgaddress{\street{2515 Speedway}, \city{Austin}, \postcode{78712}, \state{TX}, \country{US}}}

\affil[2]{\orgdiv{Oden Institute for Computational Engineering and Sciences}, \orgname{The University of Texas at Austin}, \orgaddress{\street{201 East 24th Street, Stop C0200}, \city{Austin}, \postcode{78712-1229}, \state{TX}, \country{US}}}

\abstract{The paper describes a numerical method for solving acoustic multibody
scattering problems in two and three dimensions.
The idea is to compute a highly accurate approximation to the scattering operator 
for each body through a local computation, and then use these scattering
matrices to form a global linear system.
The resulting coefficient matrix is relatively well-conditioned, 
even for problems involving a very large number of scatterers. 
The linear system is amenable to iterative solvers, and can readily be 
accelerated via fast algorithms for the matrix-vector multiplication such as the 
fast multipole method.
The key point of the work is that the local scattering matrices can be constructed
using potentially ill-conditioned techniques such as the method of fundamental 
solutions (MFS), while still maintaining scalability and numerical stability
of the global solver.
The resulting algorithm is simple, as the MFS is far simpler
to implement than
alternative techniques based on discretizing boundary integral equations using
Nystr\"om or Galerkin.}


\maketitle

\section{Introduction}
\label{sec:introduction}

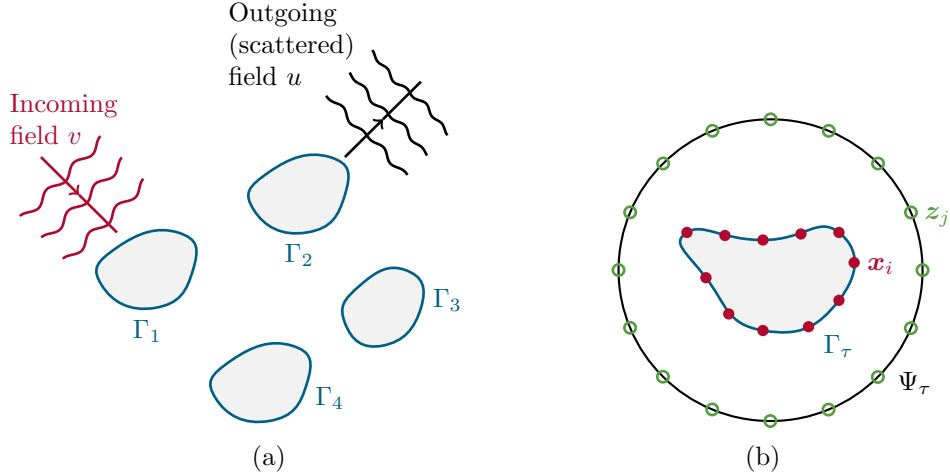
\begin{figure}
\begin{tikzpicture}
    \begin{scope}
    \filldraw[gray!10,draw=UTblue,line width=1pt] plot [smooth cycle, tension=1] coordinates {(0,0) (0.5,0.5) (0.3,1) (-0.5,0.8) (-0.7,0.3)};
    \node [UTblue] at (-0.1,-0.3) {$\Gamma_1$};
    \filldraw[gray!10,draw=UTblue,line width=1pt] plot [smooth cycle, tension=1] coordinates {(2,1) (2.5,1.5) (2.3,2) (1.5,1.8) (1.3,1.3)};
    \node [UTblue] at (1.9,0.70) {$\Gamma_2$};
    \filldraw[gray!10,draw=UTblue,line width=1pt] plot [smooth cycle, tension=1] coordinates {(3,-0.5) (3.5,0) (3.3,0.5) (2.7,0.3) (2.5,-0.2)};
    \node [UTblue] at (3.85,0.1) {$\Gamma_3$};
    \filldraw[gray!10,draw=UTblue,line width=1pt] plot [smooth cycle, tension=1] coordinates {(1.5,-1.5) (2,-1) (1.8,-0.5) (1,-0.7) (0.8,-1.2)};
    \node [UTblue] at (2.3,-1.2) {$\Gamma_4$};

    \begin{scope}[mycrimson,line width=1pt,scale=0.5,xshift=-2cm,yshift=3cm]
        \draw [->] (-1,1) -- (0,0);
        \draw [-] (0,0) -- (1,-1);
        \begin{scope}[rotate=45]
        \draw[domain=-1.5:1.5,smooth,variable=\t,samples=100] plot ({\t},{0.6-0.1*cos(deg(2*pi*\t))});
        \draw[domain=-1.5:1.5,smooth,variable=\t,samples=100] plot ({\t},{-0.2-0.1*cos(deg(2*pi*\t))});
        \draw[domain=-1.5:1.5,smooth,variable=\t,samples=100] plot ({\t},{-1-0.1*cos(deg(2*pi*\t))});
        \end{scope}
    \end{scope}
    \node [mycrimson,align=left] at (-1.2,2.5) {Incoming\\field $v$};
    
    \begin{scope}[black,line width=1pt,scale=0.5,xshift=6cm,yshift=5cm]
        \draw [->] (-1,-1) -- (0,0);
        \draw [-] (0,0) -- (1,1);
        \begin{scope}[rotate=135]
        \draw[domain=-1.5:1.5,smooth,variable=\t,samples=100] plot ({\t},{0.6-0.1*cos(deg(2*pi*\t))*sin(deg(2*\t))});
        \draw[domain=-1.5:1.5,smooth,variable=\t,samples=100] plot ({\t},{-0.2-0.1*cos(deg(2*pi*\t))*sin(deg(2.5*\t))});
        \draw[domain=-1.5:1.5,smooth,variable=\t,samples=100] plot ({\t},{-1-0.1*cos(deg(2*pi*\t))*sin(deg(3*\t))});
        \end{scope}
    \end{scope}
    \node [black,align=left] at (1.8,3.5) {Outgoing\\(scattered)\\field $u$};
    \node at (1.5,-2) {(a)};
    \end{scope}
    \begin{scope}[xshift=8cm]
    \draw[thick] (0.1,0.5) circle (2); 
    \node at (2,-1) {$\Psi_\tau$};
    \filldraw[gray!10,draw=UTblue,line width=1pt] plot [smooth cycle, tension=1] coordinates {(-1,1) (0,0.9) (1,1) (1,0.1) (0,-0.3) (-0.75,0.4)};
    \node [UTblue] at (1,-0.5) {$\Gamma_\tau$};
    \node at (0,-2) {(b)};
    \foreach \p in {(-1,1), (-0.5,0.96), (0,0.9), (0.5,0.98), (1,1), (1.2,0.6), (1,0.1), (0.6,-0.25), (0,-0.3), (-0.45,-0.08), (-0.75,0.4)} {
        \node [circle, fill=mycrimson, inner sep=1.5pt] at \p {};
    }
    \node [mycrimson] at (1.55,0.55) {$\pvct{x}_i$};
    \foreach \t in {0,0.3926,...,6.0868} {
    \node [circle, draw=UTgreen, inner sep=1.5pt, line width=1pt] at ({0.1+2*cos(deg(\t))},{0.5+2*sin(deg(\t))}) {};
    }
    \node [UTgreen] at (2.3,1.20) {$\pvct{z}_j$};
    \end{scope}
\end{tikzpicture}
\caption{
(a) The multibody scattering problem (\ref{eq:BVP}) with $T=4$ scattering
bodies,
(b) The proxy circle introduced in Section \ref{sec:geometryofSmatrix};
the scattering matrix $\mtx{S}_{\tau}$ is computed to be accurate for
interactions between $\Gamma_{\tau}$ and other objects outside of $\Psi_\tau$
}
\label{fig:geom}
\end{figure}

The numerical solution of multibody acoustic scattering problems is a recurring and important problem in scientific computing.
Given an incoming field that hits a collection of reflecting bodies
(see Figure \ref{fig:geom}(a)),
our task is to compute the resulting scattered field.
A well-known technique for solving this problem is the method of
fundamental solutions (MFS), which is often much easier to implement
than alternative methods based on either a boundary integral equation (BIE)
formulation of the problem (discretized using, say, Nystr\"om or Galerkin),
or on the direct discretization of the Helmholtz equation on an artificially
truncated domain.
The key drawback of the MFS is that it results in dense linear systems that
are invariably ill-conditioned, even to the point of being fully singular
in floating-point arithmetic.
For problems of modest scale, the ill-conditioning can be managed by
using backwards-stable least-squares solvers that rely on full matrix
factorizations such as QR or the singular value decomposition.
However, the cubic complexity of such methods greatly limits the set of
problems that can be considered.

The key contribution of the paper is a technique that allows a user to 
combine the ease of implementation and high accuracy of the MFS with the 
scalability and numerical stability of BIE-based methods for large-scale problems.

The method is based on computing an approximate ``scattering matrix'' for
each individual scatterer, using the MFS. This matrix maps an incoming field for a single
body to its reflected field. The end result is a global linear system of
the form
\begin{equation}
\hat{\vct{q}} + \mtx{S}\hat{\mtx{G}}\hat{\vct{q}} = \mtx{S}\hat{\vct{v}},
\end{equation}
where $\mtx{S}$ is a block-diagonal matrix holding the scattering matrices,
$\hat{\vct{q}}$ is a vector that holds ``equivalent sources'' that generate
the scattered field, and $\hat{\vct{v}}$ is a
representation of the externally applied incoming
field. The matrix $\mtx{S}$ is constructed via the method of fundamental
solutions using a backwards-stable solver to control a problem that is
ill-conditioned, but local. The matrix $\hat{\mtx{G}}$ can be applied rapidly
using standard fast algorithms such as the fast multipole method (FMM)
\citep{rokhlin1987,2006_rokhlin_wideband}, since its non-zero entries take the form
\begin{equation}
  \label{eq:Gmtx-entries}
\hat{\mtx{G}}(i,j) = \phi_{\kappa}(\hat{\pvct{x}}_{i} - \hat{\pvct{x}}_{j}),
\end{equation}
where $\phi_{\kappa}$ is the free-space fundamental solution of the
Helmholtz equation (cf.~\cref{eq:FundamentalSol}), and where the points 
$\{\hat{\pvct{x}}_{i}\}_{i=1}^{N}$
are collocation points on the surface of the scattering bodies.

The general notion of scattering matrices in the context of acoustic scattering is classical \citep{1954_gerjuoy_saxon,1969_waterman}, but have in the past couple of decades been used to construct fast direct solvers for boundary integral equations \citep{2013_gimbutas_greengard,2013_martinsson_smooth_BIE,2014_martinsson_simplified,2014_lai_kobayashi_greengard,2015_hao,2015_lai_kobayashi_barnett,2022_martinsson_lippmannschwinger}. For a summary, see \citep[Ch.~18]{2019_martinsson_fast_direct_solvers}. While methods based on boundary integral equations (BIE) are typically numerically stable (due to the use of second-kind integral equations) and scalable (due to the fast convergence of iterative solvers together with algorithms such as the fast multipole method), their implementation can be complicated, especially in 3D, \textit{e.g.}\ due to singular and nearly singular integrals that require special quadrature methods. For a discussion on the discretization of BIEs, see \textit{e.g.}\ \citep[Ch.~18]{2019_martinsson_fast_direct_solvers}, \citep{2013_martinsson_nystrom_final}; special quadrature methods for BIEs are an active area of research \citep{2013_klockner_qbx,2020_wala_qbx,2022_zhu_veerapaneni,2025_kaneko_duraiswami}.

The method of fundamental solutions (MFS) is an alternative method that is easy to implement and requires no special quadrature. The basic ideas of the formulation of the MFS were first proposed by Kupradze and Aleksidze \cite{kupradze1964method,1964_kupradze_b}. The method has also been called the ``charge simulation method'' \citep{1988_katsurada}, the ``source simulation technique'' \citep{1992_mohsen} and the ``method of auxiliary sources'' \citep{2002_shubitidze}. For an overview of the MFS, see \textit{e.g.}\ \cite{2003_fairweather_mfs,cheng2020overview}. Important analysis of the MFS applied to acoustic scattering was provided by Barnett and Betcke \cite{2008_barnett_betcke}. The method was studied extensively in this context by Liu and Barnett \cite{2016_liu_dissertation,2016_liu_barnett}, including for domains with corners.
The ``lightning solver'' proposed by Gopal and Trefethen \cite{2019_gopal_trefethen_lightning_sinum,2019_gopal_helmholtz} is in many ways similar to the MFS, but uses a different set of basis functions. The scattering matrices, which we here compute using the MFS, could equally well be computed using a lightning solver.

As mentioned, the ill-conditioning of the MFS limits the size of the systems it can be applied to. Recent work has striven to lift this restriction. For instance, in work by Antunes et al.\ \cite{2022_antunes,2024_antunes_helmholtz_2d,2025_santos_preprint}, the MFS basis functions are expanded and orthogonalized in order to define a more well-conditioned basis. However, this technique has not yet been applied to exterior Helmholtz problems, which is the topic of interest of this paper. An entirely different approach is taken by Stein and Barnett \cite{2022_stein_barnett}, where the MFS is used within a BIE method to remove the need for traditional special quadrature. With a block-diagonal preconditioner, a 2D Helmholtz problem with 1000 scatterers and around 76 wavelengths across the whole geometry (around $1/3$ wavelength across one scatterer) can be solved to 10 digits of accuracy in 1331 GMRES iterations. Unlike our proposed method, this method has no compression of the global linear system.

For the Stokes problem, the block-diagonal preconditioner strategy from \cite{2016_liu_barnett,2022_stein_barnett} has been applied to turn the global linear system associated with the MFS into a square and more well-conditioned system \citep{2025_broms_close_spheres,2025_broms_large_scale_3d,2024_broms_dissertation}. In particular, Broms et al.\ \cite{2025_broms_large_scale_3d} consider a 3D Stokes problem with 10,000 ellipsoids, which is solved to 5 digits of accuracy in only 7 GMRES iterations. (Stokes problems are inherently more amenable to GMRES than the Helmholtz problem considered above.) Jordan and Lockerby \cite{2025_jordan_lockerby} apply the MFS to single particles in Stokes flow, together with an iterative scheme based on block Gauss-Seidel to resolve the particle-particle interactions. Computations with up to 6000 ellipsoids are considered, with long-range interactions replaced by a single point source per particle. The accuracy of their method, which scales quadratically in the number of particles, is however not studied carefully in the multi-particle case.

As far as we are aware, the MFS has not previously been used to compute scattering matrices in the context of acoustic scattering, and this is our major contribution in this paper. Furthermore, we demonstrate with numerical examples that the resulting linear system is reasonably well-conditioned (no worse than a well-conditioned BIE formulation) and can be solved efficiently using iterative methods.

The paper is organized as follows: 
\cref{sec:overview} provides a high-level overview of the proposed solver. 
\cref{sec:MFS} reviews some preliminaries necessary to our solver, namely the method of fundamental solutions (MFS). 
\cref{sec:constructS} presents the linear algebraic formulation of building the local scattering matrices. 
\cref{sec:NumericalResults} includes our numerical experiments in both two and three dimensions, and the final section summarizes the main findings.

\section{Overview of the proposed solver}
\label{sec:overview}

This section provides a high-level description of the proposed method.
Implementation details will follow in 
Sections~\ref{sec:MFS} and \ref{sec:constructS}.

\subsection{Problem formulation}
Let $\Omega$ be the domain external to a boundary $\Gamma$ that consists
of a union $\Gamma = \bigcup_{\tau=1}^{T}\Gamma_{\tau}$ of $T$ disjoint scattering surfaces, as shown in Figure \ref{fig:geom}(a).
We then seek to solve the Helmholtz problem
\begin{equation}
\label{eq:BVP}
\left\{\begin{aligned}
-\Delta u(\pvct{x}) - \kappa^{2}\,u(\pvct{x}) =&\ 0 \qquad &\pvct{x} \in \Omega,\\
u(\pvct{x}) =&\ v(\pvct{x}) \qquad &\pvct{x} \in \Gamma,\\
\lim_{|\pvct{x}| \rightarrow\infty}|\pvct{x}|\bigl(\partial_{|\pvct{x}|}u(\pvct{x}) - i\kappa\,u(\pvct{x})\bigr) =&\ 0,
\end{aligned}\right.
\end{equation}
where $v$ is a \textit{global incoming field} that itself satisfies the free-space Helmholtz equation, \textit{i.e.}, $-\Delta v(\pvct{x}) - \kappa^{2}\,v(\pvct{x}) = 0$, in some neighborhood of $\Gamma$, and $u$ is the scattered (outgoing) field.%
\footnote{More specifically, $v(\pxx) = -u_\text{inc}(\pxx)$ for an incoming field $u_\text{inc}$, and the boundary condition comes from the condition that the total field $u_\text{tot}(\pxx) = u(\pxx) + u_\text{inc}(\pxx) = 0$, $\pxx \in \Gamma$, for sound-soft objects.}
We focus on Dirichlet boundary data (\textit{i.e.}, sound-soft scattering objects) for concreteness, but the method can handle
other boundary conditions with no additional difficulty.

\subsection{Local scattering matrices and forming a global system}
\label{sec:introS}
The solver we describe is based on the notion of \textit{scattering matrices}.
To describe this concept, let us consider a single scatterer $\Gamma_{\tau}$,
as pictured in Figure \ref{fig:geom}(b).
Let $\hat{\vct{w}}_{\tau}$ denote some discrete representation of the
\textit{local incoming field} on $\tau$, which is defined as the sum of the externally
applied field $v$, and the scattered fields from the other $T-1$ scatterers.
In response to the incoming field, the scatterer $\tau$ generates an
\textit{outgoing field} $u_{\tau}$ that we represent using some
vector $\hat{\vct{q}}_{\tau}$.
In the physics literature, and previous work such as \cite{2014_lai_kobayashi_greengard,2015_lai_kobayashi_barnett}, it is standard practice to have the vectors
$\hat{\vct{w}}_{\tau}$ and $\hat{\vct{q}}_{\tau}$
hold the  expansion coefficients (Fourier in 2D, and spherical harmonics in 3D)
of the incoming and outgoing fields, respectively.
In contrast, we will use collocation at some set of points on the scattering surface to
identify the incoming field, and will then represent the outgoing field using
``equivalent sources'' at the same set of points on $\Gamma_{\tau}$.

The \textit{scattering matrix} for $\tau$ is then defined as the linear map
\begin{equation}
\label{eq:defS}
\hat{\vct{q}}_{\tau} = \mtx{S}_{\tau}\hat{\vct{w}}_{\tau}.
\end{equation}
If all matrices $\{\mtx{S}_{\tau}\}_{\tau=1}^{T}$ are available, it
is easy to build a global system for determining all the vectors
$\{\hat{\vct{q}}_{\tau}\}_{\tau=1}^{T}$. Simply observe that for any
scatterer $\tau$, its incoming field takes the form
\begin{equation}
\label{eq:w_in}
\hat{\vct{w}}_{\tau} = \hat{\vct{v}}_{\tau} - \sum_{\sigma\neq \tau}\hat{\mtx{G}}_{\tau,\sigma}\hat{\vct{q}}_{\sigma},
\end{equation}
where $\hat{\mtx{G}}_{\tau,\sigma}$ is the matrix whose entries
are the fundamental solution evaluated between points on
$\Gamma_\sigma$ and $\Gamma_\tau$; it converts an outgoing expansion
for $\Gamma_{\sigma}$ to an incoming expansion on $\Gamma_{\tau}$.
Combining \cref{eq:defS,eq:w_in}, we get
\begin{equation}
\label{eq:mbody1}
\hat{\vct{q}}_{\tau} + \mtx{S}_{\tau}\sum_{\sigma \neq \tau}\hat{\mtx{G}}_{\tau,\sigma}\hat{\vct{q}}_{\sigma} = \mtx{S}_{\tau}\hat{\vct{v}}_{\tau},
\qquad \tau \in \{1,\,2,\,\dots,\,T\}, 
\end{equation}
which we write compactly as
\begin{equation}
\label{eq:mbody2}
\bigl(\mtx{I} + \mtx{S}\hat{\mtx{G}}\bigr)\hat{\vct{q}} = \mtx{S}\hat{\vct{v}},
\end{equation}
where $\mtx{S}$ is a block-diagonal matrix holding the scattering matrices,
and $\hat{\mtx{G}}$ is a dense matrix holding the off-diagonal blocks $\hat{\mtx{G}}_{\tau,\sigma}$.

The formulation (\ref{eq:mbody2}) of the multibody scattering problem has several advantages:
\begin{enumerate}
\item The linear system (\ref{eq:mbody2}) tends to be well-conditioned even for
problems involving large numbers of scatterers, and is well suited for iterative
methods. Numerical examples in Sections \ref{sec:2d-results} and \ref{sec:3d-results}
demonstrate that the local scatterers may even contain cavities.
\item Iterative solvers for (\ref{eq:mbody2}) are computationally efficient,
since the matrix $\bigl(\mtx{I} + \mtx{S}\hat{\mtx{G}}\bigr)$
can be applied to vectors efficiently. The matrix $\mtx{S}$ is block-diagonal, and the matrix $\hat{\mtx{G}}$ can be applied using fast algorithms such as
the fast multipole method \citep{rokhlin1987,2006_rokhlin_wideband}.
\item The number of degrees of freedom in the global linear system (\ref{eq:mbody2})
is determined by the numerical ranks of interaction between the scatterers, 
\textit{independent of the internal geometric complexity of each
individual scatterer}. This is in contrast to regular BIE-based methods, where the
size of the global system is determined by the need to resolve
all local features. In our approach, we still need
to introduce enough variables to fully resolve the problem locally, but this
computation involves only temporary variables, and it can be executed 
efficiently and accurately using dense linear algebra.
\end{enumerate}

\subsection{Computing the local scattering matrices}
The technique for formulating a multibody scattering problem through local
scattering matrices was explored by, \textit{e.g.}, \cite{2013_martinsson_smooth_BIE,2022_martinsson_lippmannschwinger,2014_martinsson_simplified}.
In these works, the starting point was a BIE formulation for the local scattering
problems. For a single scatterer $\Gamma_{\tau}$ that requires $m$ degrees of
freedom to resolve%
\footnote{For simplicity, we write $m$ and $k$ throughout this paper, but of course, these numbers could be different for each scatterer $\Gamma_\tau$, \textit{i.e.}, $m = m_\tau$ and $k = k_\tau$.}, the matrix $\mtx{S}_{\tau}$ then takes the form
\begin{equation}
\label{eq:localold}
\begin{array}{ccccccccccccccc}
\mtx{S}_{\tau} &=& \mtx{V}_{\tau}^{*}&\mtx{A}_{\tau}^{-1}&\mtx{U}_{\tau},\\
k \times k && k\times m & m\times m & m\times k
\end{array}
\end{equation}
where $k$ is the numerical rank of interaction between $\Gamma_{\tau}$ and the
rest of the geometry. The matrix $\mtx{A}_{\tau}$ is the discretization of a BIE
on $\Gamma_{\tau}$. The matrices $\mtx{U}_{\tau}$ and $\mtx{V}_{\tau}$ are well-conditioned
matrices that translate between the compact representations of the incoming and outgoing
fields, and the fully resolved local problem.

The key novelty of the present paper is the introduction of a new technique for
forming the scattering matrix $\mtx{S}_{\tau}$ using the method of fundamental
solutions (MFS). The advantage of this method is that it is very simple to implement,
as there is no need to introduce specialized quadrature rules for singular integrals.
The drawback is that the local problem becomes ill-conditioned, and requires expensive
dense linear algebraic techniques. To be precise, the formula (\ref{eq:localold})
gets replaced by
\begin{equation}
\label{eq:localnew}
\begin{array}{ccccccccccccccc}
\mtx{S}_{\tau} &=& \mtx{V}_{\tau}^{*}&\mtx{A}_{\tau}^{\dagger}&\mtx{U}_{\tau},\\
k \times k && k\times n & n\times m & m\times k
\end{array}
\end{equation}
where $\mtx{A}_{\tau}$ is now a rectangular matrix of size $m\times n$ (with $m > n$), 
and where $(\cdot)^{\dagger}$ denotes a pseudo-inverse,
which in practice could be applied using either a truncated
singular value decomposition or a QR factorization.

\subsection{Solving the global system}

To reiterate, the global system is in the form of (\ref{eq:mbody2}):
\begin{equation*}
\bigl(\mtx{I} + \mtx{S}\hat{\mtx{G}}\bigr)\hat{\vct{q}} = \mtx{S}\hat{\vct{v}}.
\end{equation*}
Empirical evidence indicates that the system is reasonably well-conditioned and can be solved using an iterative solver such as GMRES \citep{1986_saad_gmres}, which will converge in a small number of iterations relative to the system size.
Detailed numerical experiments are described in \cref{sec:2d-results}, Example~\ref{ex:2d-results-scaling}, and \cref{sec:3d-results}, Example~\ref{ex:scaling-3D}.
The solving process can be accelerated via fast algorithms for the matrix-vector multiplication, such as the fast multipole method (FMM).

Once the equivalent sources $\hat{\vct{q}}$ of the scattered field have been determined, the scattered field $u(\pvct{x})$ itself can be evaluated for any $\pvct{x} \in \Omega$ using the standard fundamental solution of the Helmholtz equation (see \cref{sec:geometryofSmatrix}). Again, this evaluation can be accelerated using fast methods such as the FMM. If the evaluation point $\pvct{x}$ lies within the proxy surface $\Phi_\tau$ of one (or more) of the scatterers, the contribution from that scatterer $\Gamma_\tau$ is computed by reconstructing the vector of full sources $\vct{q}$ for $\Gamma_\tau$, which is a local computation. This reconstruction is described in \cref{sec:reconstruction}.

\section{Method of fundamental solutions}
\label{sec:MFS}

In this section, we briefly review how the classical
``method of fundamental solutions (MFS)'' can be applied
to solve a \emph{single-body} exterior scattering problem.

\subsection{Problem formulation}
Let us consider the single-body scattering problem, cf.~\cref{eq:BVP},
\begin{equation}
\label{eq:BVPlocal}
\left\{\begin{aligned}
-\Delta u(\pvct{x}) - \kappa^{2}\,u(\pvct{x}) =&\ 0 \qquad &\pvct{x} \in \Omega,\\
u(\pvct{x}) =&\ w(\pvct{x}) \qquad &\pvct{x} \in \Gamma,\\
\lim_{|\pvct{x}| \rightarrow\infty}|\pvct{x}|\bigl(\partial_{|\pvct{x}|}u(\pvct{x}) - i\kappa\,u(\pvct{x})\bigr) =&\ 0,
\end{aligned}\right.
\end{equation}
for a domain $\Omega$ exterior to a single scatterer $\Gamma = \partial\Omega$, where $w$ is the (local) incoming field.
The MFS formulation relies crucially on the 
fundamental solution of the Helmholtz operator $-\Delta - \kappa^2 I$, which is given by 
\begin{equation}
\label{eq:FundamentalSol}
    \phi_{\kappa} (\pxx) = \left\{\begin{aligned}
    & \frac{i}{4} H_0^{(1)} (\kappa |\pxx|) \quad& \text{in 2D}, \\[5pt]
    & \frac{1}{4\pi} \frac{e^{i\kappa |\pxx|}}{|\pxx|} \quad& \text{in 3D},
\end{aligned}\right.
\end{equation}
where $H_0^{(1)}$ is the 0th-order Hankel function of the first kind.

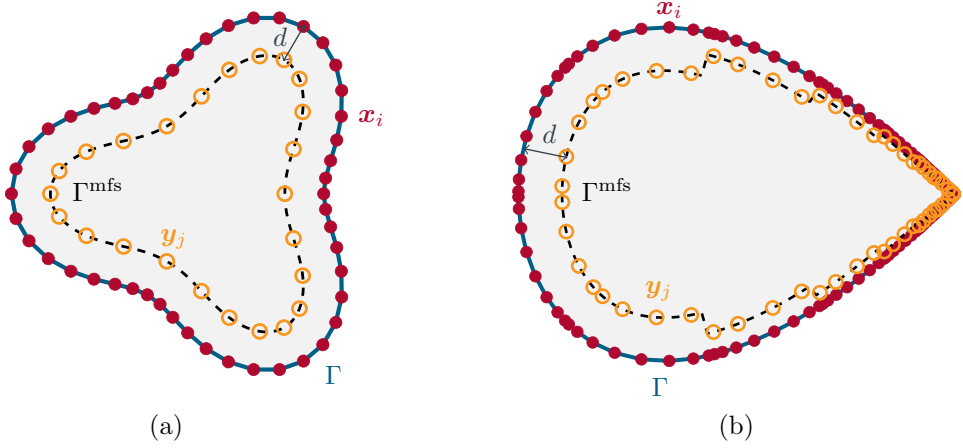
\begin{figure}
\begin{tikzpicture}
    \begin{scope}[scale = 2]
    \node [inner sep=0pt] at (0,0) {\includegraphics[width=5cm,trim=6cm 9cm 6cm 9cm,clip]{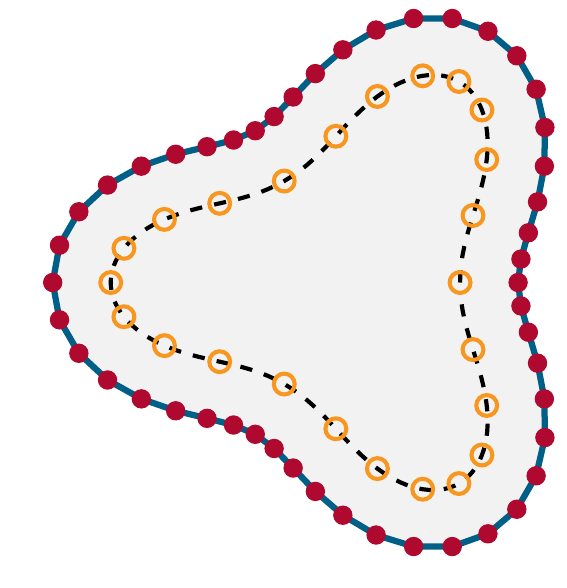}};
    \node [UTblue] at (1.1,-1.15) {$\Gamma$};
    \node at (-0.45,0.08) {$\Gamma^\text{mfs}$};
    \node [mycrimson] at (1.35,0.55) {$\pvct{x}_i$};
    \node [UTorange] at (0.05,-0.25) {$\pvct{y}_j$};
    \draw [UTdark, line width=0.5pt] [<->] (0.91,1.17) -- (0.77,0.93);
    \node [UTdark] at (0.75,1.09) {$d$};
    \node at (0,-1.5) {(a)};
    \end{scope}

    \begin{scope}[scale = 2, xshift=3.75cm]
    \node [inner sep=0pt] at (0,0) {\includegraphics[width=6cm,trim=5cm 9cm 4.5cm 9cm,clip]{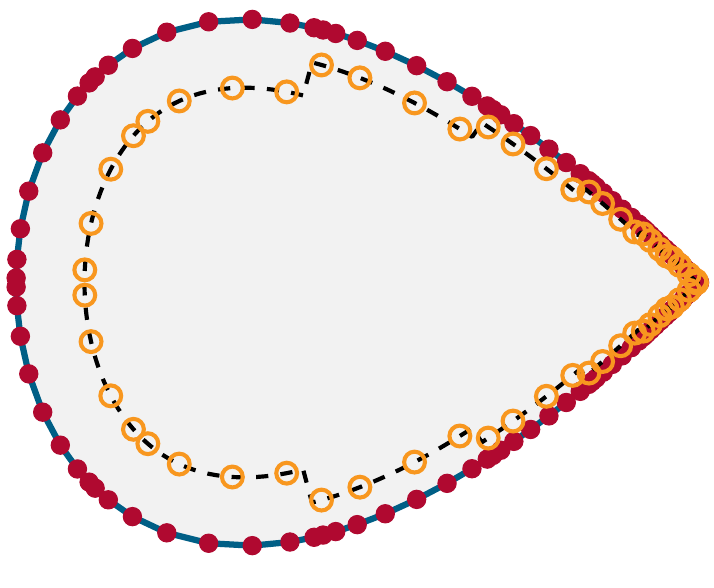}};
    \node [UTblue] at (-0.5,-1.23) {$\Gamma$};
    \node at (-0.85,0.08) {$\Gamma^\text{mfs}$};
    \node [mycrimson] at (-0.44,1.27) {$\pvct{x}_i$};
    \node [UTorange] at (-0.51,-0.59) {$\pvct{y}_j$};
    \draw [UTdark, line width=0.5pt] [<->] (-1.11,0.29) -- (-1.4,0.35);
    \node [UTdark] at (-1.23,0.45) {$d$};
    \node at (0,-1.5) {(b)};
    \end{scope}
\end{tikzpicture}
\caption{A single-body scattering problem for (a) a smooth object or (b) an object with a corner}
\label{fig:MFS}
\end{figure}

In the MFS, the general idea is to pick a set of points $\{\pyy_{j}\}_{j=1}^{n}$ \textit{inside}
the domain enclosed by the boundary $\Gamma$, \textit{i.e.}, in $\mathbb{R}^d \setminus (\Omega \cup \Gamma)$, as shown in Figure~\ref{fig:MFS}.
We refer to these points as the \textit{MFS sources}, and place them on an MFS surface (curve) $\Gamma^\text{mfs}$ a distance $d$ from the boundary $\Gamma$.
We then seek a solution $u$ to (\ref{eq:BVPlocal}) of the form
\begin{equation}
\label{eq:MFSansatz}
u(\pxx) \approx \sum_{j=1}^{n}\phi_{\kappa}(\pxx - \pyy_{j})\,q_{j},\qquad\pxx \in \Omega,
\end{equation}
where $\{q_{j}\}_{j=1}^{n}$ is a set of ``MFS source strengths'' that we will need to determine.
Observe that any function $u$ of the form (\ref{eq:MFSansatz}) will
satisfy both the partial differential equation in $\Omega$, and the radiation condition, so
all that remains is to attempt to satisfy the boundary condition
in some approximate sense. To achieve this, we pick a set of collocation
points $\{\pxx_{i}\}_{i=1}^{m} \subset \Gamma$, and then attempt to
enforce that
\begin{equation}
\label{eq:MFSlinsys}
w(\pxx_{i}) = \sum_{j=1}^{n}\phi_{\kappa}(\pxx_{i} - \pyy_{j})\,q_{j},
\qquad i \in \{1,\,2,\dots,\,m\}.
\end{equation}
We write \cref{eq:MFSlinsys} compactly as
\begin{equation}
\label{eq:MFSlinsy2}
\begin{array}{ccccccccc}
\mtx{A} & \vct{q} &=& \vct{w},\\
m\times n & n\times 1 && m\times 1
\end{array}
\end{equation}
where $\mtx{A}$ is the $m\times n$ matrix with entries
\begin{equation}
\mtx{A}(i,j) = \phi_{\kappa}(\pxx_{i} - \pyy_{j}).
\end{equation}

The MFS almost always results in a linear system that is either highly
ill-conditioned, or fully singular. 
This ill-conditioning is a result of the non-uniqueness of the set of MFS sources. 
While the ill-conditioning poses challenges, these can be overcome by using 
backward-stable solvers for the least-square problem such as a truncated 
singular value decomposition, or a QR factorization.
These techniques work well, but require $O(mn\,\min(m,n))$ flops, which
limits the classical MFS to small and medium scale problems.
Here, we follow standard practice of using more collocation points than
source points (so that $m > n$), which tends to improve both stability
and accuracy \citep{2016_liu_dissertation,cheng2020overview}.
What is important is to make sure that MFS source locations are chosen
in such a way that the linear system  (\ref{eq:MFSlinsy2}) admits a solution
with small residual (meaning that the Ansatz (\ref{eq:MFSansatz}) is well
capable of resolving the solution that is sought).

The accuracy of the MFS solution is highly sensitive to the placement of the MFS sources. According to \cite{2008_barnett_betcke}, the source points should be placed such that the analytical continuation of the solution into the scatterer has no singularities up until the MFS surface $\Gamma^\text{mfs}$. Otherwise, the MFS source strengths become large, which leads to a loss of digits. As long as there are no singularities between $\Gamma$ and $\Gamma^\text{mfs}$, a larger distance between the two tends to give faster convergence as $n \to \infty$.

Since the locations of singularities of the solution are in general not known, the placement of the MFS surface is often based on heuristics. Finding suitable source locations is crucial, since badly selected sources may lead to a significant increase in the error, often by several orders of magnitude \cite{2009_alves,li2013mfs,chen2016choosing,2016_liu_dissertation}. Below, we outline how the MFS sources are selected in the different cases.

\subsection{Selection of points for smooth contours in 2D}
\label{sec:mfs-smooth-2d}

In 2D, experiments by \cite{li2013mfs} showed that a uniform distribution of collocation points around the boundary $\Gamma$, and a uniform distribution of the MFS source points along a curve at a fixed distance from $\Gamma$, produces the most accurate and stable results. We follow this suggestion, but, for simplicity, we place both sets of points equidistantly in parameter space rather than in arclength. This way of selecting the points was considered also by \cite{2016_liu_dissertation}.

Using experiments, we determine a suitable ratio between the number of collocation points $m$ and MFS sources $n$. Setting the ratio to 2 leads to a good tradeoff between stability and cost, so we fix this value for all the numerical experiments, \textit{i.e.}, $ n = m/2 $. At the locations where MFS sources should be placed, we compute the normal vectors to the physical boundary, pointing away from $\Omega$ (\textit{i.e.}, into the scatterer). MFS sources are then placed in the direction of the normal vectors at a distance $d$ from the contour.

We determine the MFS distance $d$ using numerical experiments.
As an alternative, the Leave-One-Out Cross Validation (LOOCV) method could be adopted to assist in the search for an optimal $d$, as introduced by \cite{chen2016choosing,cheng2020overview}. However, we do not use this method here, since our focus is not on finding optimal MFS parameters.

(We also tried using Gauss--Legendre panels instead of equidistant points for smooth contours. While both options work, we found that equidistant points gave slightly lower errors for the same total number of points.)

\subsection{Selection of points for non-smooth contours in 2D}
\label{sec:mfs-nonsmooth-2d}

Domains with corners are usually challenging to deal with, as they will lead to solutions with a singularity or sharp gradients at the corner. To capture the singular behavior more accurately, a high density of nodes is needed at the corner. Therefore, adaptive mesh refinements are often used for discretization in the vicinity of corners. Due to its adaptivity, the MFS gives reasonably accurate results when applied to problems involving boundary singularities.

For non-smooth contours, we use Gauss--Legendre panels for discretization. Dyadic refinement of both collocation points and MFS sources is needed close to geometrical singularities such as corners, as shown in Figure \ref{fig:MFS}(b). The MFS distance is adjusted locally based on the dyadic refinement, as seen in the figure. This way of selecting the points is essentially identical to the ``method of refined panels'' of \cite[Sec.~5.1]{2016_liu_dissertation}.

Our treatment of corner singularities also has some similarities to the ``lightning solver'' introduced by \cite{2019_gopal_trefethen_lightning_sinum,2019_gopal_helmholtz}. Similar to the MFS, the lightning solver approximates solutions via finite sums, and deals with corners by exponentially clustering source points at the corner. It differs from the MFS in that the smooth part of the solution is represented by a polynomial expansion, while the corner singularity is represented by extra dipole terms along a line emanating from the corner, rather than monopoles along the boundary as in the MFS.

\subsection{Selection of points for surfaces in 3D}

In the three-dimensional case, we consider only smooth surfaces, and place both collocation points and MFS sources equidistantly in each parameter, forming a tensor product grid in parameter space. The ratio between the number of collocation points and MFS sources is $2$ along each parameter, which means that for the surface as a whole, the ratio is $4$, \textit{i.e.}, $n=m/4$. The MFS distance $d$ is again fixed for the surface as a whole, and selected based on experiments.

Since the MFS is based on a direct linear combination of fundamental solutions, it is in many ways simpler than alternative techniques based on discretizing boundary integral equations. Such approaches require defining an underlying quadrature rule for the boundary, and also require special quadrature for singular and near-singular integration, all of which can be avoided with the MFS. This advantage makes the MFS mathematically and computationally simpler, especially for surfaces in 3D.

\section{Construction of local scattering matrices}
\label{sec:constructS}

In this section, the focus remains on the single-body scattering
problem (\ref{eq:BVPlocal}). We now describe how to build a
scattering matrix $\mtx{S}_{\tau}$ that maps an incoming field
to the scattered field, accurate to some requested precision $\varepsilon$.
The cost of the method we describe is $O(m^3)$ where $m$ is the number of 
degrees of freedom required to fully resolve the local problem.

\subsection{Problem formulation}
\label{sec:geometryofSmatrix}
In order to numerically construct an approximate scattering matrix,
we must specify a requested precision $\varepsilon$, as well as
what fields it should be able to resolve. For the latter,
we place what we call a \textit{proxy surface} $\Psi_\tau$ around $\Gamma_{\tau}$,
cf.~Figure \ref{fig:geom}(b).
Then our objective is to handle any incoming field $w$ that satisfies
the free space Helmholtz equation \textit{inside $\Psi_\tau$}. Additionally, our
construction should allow for the evaluation of the scattered (outgoing)
field (to within
precision $\varepsilon$) at any point \textit{on or outside $\Psi_\tau$}.
In other words, the scattering matrix should be able to handle all
interactions with scatterers that are located outside of $\Psi_{\tau}$.

To represent the incoming field, we will select a set
of ``skeleton'' points $\{\hpxx_{i}\}_{i=1}^{k} \subset \Gamma_{\tau}$
that form good interpolation points. In other words, given the values
of any incoming field $w$ at the points $\{\hpxx_{i}\}_{i=1}^{k}$, 
we can stably reconstruct $w$ (to precision $\varepsilon$) at any other
point on $\Gamma_{\tau}$.

Due to the symmetry of the problem, the skeleton points
$\{\hpxx_{i}\}_{i=1}^{k} \subset \Gamma_{\tau}$ chosen for the incoming
field will necessarily also work well to represent any outgoing field.
In other words, given any outgoing field $u_\tau$, we can find a set of
equivalent sources $\{\hat{q}_{i}\}_{i=1}^{k}$ such that
\begin{equation}
\label{eq:out_skel}
u_{\tau}(\pxx) \stackrel{\varepsilon}{=} \sum_{i=1}^{k}\phi_{\kappa}(\pxx - \hpxx_{i})\,\hat{q}_{i}
\end{equation}
for any point $\pxx$ on or outside of $\Psi_\tau$.

\subsection{Determining skeleton points and interpolation operators}
\label{sec:building-interpolation}

We next describe how to determine a good set of skeleton points, and the associated
interpolation operators (the matrices $\mtx{U}_{\tau}$ and $\mtx{V}_{\tau}$ in 
(\ref{eq:localnew})). 
Our starting point is a local surface $\Gamma_{\tau}$ with an associated
set of discretization nodes $\{\pxx_{i}\}_{i=1}^{m}$.
We are also given a proxy surface $\Psi_{\tau}$.

The first step is to place points $\{\pzz_{i}\}_{i=1}^{p}$ uniformly on the proxy surface $\Psi_{\tau}$. We then build an outgoing $p \times m$ matrix $\mtx{B}_{\tau}$ with entries
\begin{equation}
\label{eq:B-matrix}
\mtx{B}_{\tau}(i,j) = \phi_{\kappa}(\pzz_{i} - \pxx_{j}).
\end{equation}
Multiplying the matrix $\mtx{B}_{\tau}$ by a vector computes an outgoing field on the proxy surface $\Psi_\tau$ generated by sources on $\Gamma_\tau$.
The key step is now to perform a (rank-revealing) column-pivoted QR factorization 
on the matrix $\mtx{B}_{\tau}$, stopping once a precision of $\varepsilon$ has been
attained.
Let $k$ denote the rank of the factorization when it is halted, and let 
$\hat{I}_{\tau} \subseteq \{1,\,2,\,\dots,\,m\}$ denote an index vector
that marks the chosen columns. 
Then $\{\pxx_{i}\}_{i \in \hat{I}_{\tau}}$
is the collection of skeleton points on $\Gamma_{\tau}$ (for the outgoing field).

To build the interpolation matrix, we convert the partial QR factorization
of $\mtx{B}_{\tau}$ into an \textit{interpolative decomposition (ID)}, as
described in \cite[Ch.~3]{2019_martinsson_fast_direct_solvers}. This results
in a factorization
\begin{equation}
\label{eq:B_fac}
\begin{array}{ccccc}
  \mtx{B}_{\tau} & \stackrel{\varepsilon}{=} & \mtx{B}_{\tau}(\colon,\hat{I}_\tau) 
  & \mtx{Z}_\tau^*,\\
  m\times p && m\times k & k\times p
\end{array}
\end{equation}
where $\mtx{Z}_{\tau}$ is a $p\times k$ ``interpolation matrix'' that contains
the $k\times k$ identity matrix as a submatrix.

Next, let us consider the incoming skeleton. We build the $p\times m$ matrix
$\mtx{E}_{\tau} = \mtx{B}_{\tau}^{\mathsf T}$ with entries
\begin{equation}
\mtx{E}_{\tau}(i,j) = \phi_{\kappa}(\pxx_{i} - \pzz_{j}).
\end{equation}
This matrix maps sources on $\Psi_{\tau}$ to a potential on $\Gamma_{\tau}$.
We observe that any incoming field is then (to precision $\varepsilon$)
inside the range of $\mtx{E}_{\tau}$. Since $\mtx{E}_{\tau}$ is the 
transpose of $\mtx{B}_{\tau}$ (this follows from the symmetry of the fundamental solution $\phi_\kappa$), we immediately get a decomposition of it as
\begin{equation}
    \mtx{E}_{\tau} = \mtx{B}_{\tau}^\mathsf{T}  \stackrel{\varepsilon}{=} (\mtx{Z}_\tau^*)^\mathsf{T} (\mtx{B}_{\tau}(\colon,\hat{I}_\tau))^\mathsf{T}
     = \overline{\mtx{Z}}_{\tau} \mtx{E}_{\tau}(\hat{I}_\tau,\colon) 
     = \mtx{U}_{\tau} \mtx{E}_{\tau}(\hat{I}_\tau,\colon),
\end{equation}
where $\mtx{U}_{\tau} = \overline{\mtx{Z}}_{\tau}$ is the complex conjugate of $\mtx{Z}_{\tau}$. This shows that the same skeleton points can be used for the incoming and outgoing fields.

\subsection{Solving the local equilibrium problem}

Introduce the $m\times n$ matrix $\mtx{A}_{\tau}$ from Section~\ref{sec:MFS} that encodes the local MFS problem \eqref{eq:BVPlocal}:
\begin{equation}
\mtx{A}_{\tau}(i,j) = \phi_{\kappa}(\pvct{x}_{i} - \pvct{y}_{j}),
\end{equation}
where $\{\pvct{y}_{j}\}$ is the set of MFS source locations.
Note that $u(\pvct{x}) = w(\pvct{x})$ on $\Gamma_\tau$ in this local problem, where $w$ is the incoming field and $u$ is the outgoing field.
The local MFS system is then
\begin{equation}
\label{eq:local_mfs}
    \mtx{A}_{\tau} \vct{q}_{\tau}^{\text{mfs}} = \vct{w}_{\tau},
\end{equation}
where $\vct{q}_{\tau}^{\text{mfs}}$ contains the MFS source strengths for $\Gamma_\tau$ and $\vct{w}_\tau$ is the incoming field evaluated at the collocation points of $\Gamma_\tau$.

Let $\vct{q}_{\tau}^{\text{proxy}}$ represent the charges at the proxy points that would generate the field $w_\tau$ on $\Gamma_{\tau}$. Considering the incoming field on $\Gamma_{\tau}$,
\begin{equation}
    \vct{w}_{\tau} =  \mtx{E}_{\tau} \vct{q}_{\tau}^{\text{proxy}} \stackrel{\varepsilon}{=} \mtx{U}_{\tau}  \mtx{E}_{\tau}(\hat{I}_{\tau},\colon) \vct{q}_{\tau}^{\text{proxy}} =  \mtx{U}_{\tau}\hat{\vct{w}}_{\tau}.
\end{equation}
This shows that
\begin{equation}
\label{eq:inc_field}
    \vct{w}_{\tau} \stackrel{\varepsilon}{=} \mtx{U}_{\tau}\hat{\vct{w}}_{\tau},
\end{equation}
where $\hat{\vct{w}}_{\tau} = \vct{w}_\tau(\hat{I}_\tau)$ is the incoming field evaluated at the skeleton points, as introduced in Section~\ref{sec:overview}. Thus, the matrix $\mtx{U}_\tau$ lets us reconstruct the full incoming field $\vct{w}_\tau$ given the compressed representation $\hat{\vct{w}}_\tau$.

Combining \cref{eq:local_mfs,eq:inc_field}, we get the equation
\begin{equation}
\label{eq:mfs-source-solution}
\vct{q}_{\tau}^{\text{mfs}} \approx \mtx{A}_{\tau}^{\dagger}\mtx{U}_{\tau}\hat{\vct{w}}_{\tau}
\end{equation}
for computing the MFS sources, where $\mtx{A}_{\tau}^{\dagger}$ denotes the pseudo-inverse of $\mtx{A}_{\tau}$.

\subsection{Building the translation operator}
The solution \eqref{eq:mfs-source-solution} to the local equilibrium problem above is defined at the MFS source points inside the region enclosed by $\Gamma_\tau$, but we would like to define the equivalent sources $\hat{\vct{q}}$ on a subset of the grid points on $\Gamma_\tau$, as mentioned in \cref{sec:building-interpolation} above. (The equivalent sources \emph{could} be defined in a subset of the MFS points, but that leads to a less well-conditioned global equilibrium equation, as discussed in \cref{sec:build-scat} below.) Thus, we need a way to map MFS sources to sources on $\Gamma_\tau$. We use the proxy surface $\Psi_\tau$ to do this, as follows. (The idea is similar to the off-surface check surface used in the QFS-D method of \cite{2022_stein_barnett}.)

Let us introduce the $p\times n$ matrix $\mtx{D}_{\tau}$ with entries
\begin{equation}
\mtx{D}_{\tau}(i,j) = \phi_{\kappa}(\pzz_{i}-\pyy_{j}).
\end{equation}
This matrix represents a field generated by the MFS sources evaluated on the proxy surface. The $m\times n$ translation matrix $\mtx{C}_{\tau}$ is constructed
as the (stabilized) solution to the least squares problem
\begin{equation}
\begin{array}{cccccccccccccccc}
\mtx{B}_{\tau} & \mtx{C}_{\tau} &=& \mtx{D}_{\tau}, \\
p\times m & m\times n && p\times n
\end{array}
\end{equation}
where $\mtx{B}_\tau$ is the matrix from \cref{eq:B-matrix}.
Considering the outgoing field on the proxy surface $\Psi_{\tau}$,
\begin{equation}
\begin{array}{cccccccccccccccc}
\mtx{B}_{\tau} \mtx{C}_{\tau} \vct{q}_{\tau}^{\text{mfs}} = \mtx{D}_{\tau} \vct{q}_{\tau}^{\text{mfs}} =: \mtx{B}_{\tau} \vct{q}_{\tau},
\end{array}
\end{equation}
so we should set $\vct{q}_\tau = \mtx{C}_{\tau} \vct{q}^\text{mfs}_\tau$, where $\vct{q}_\tau$ is the vector of sources on $\Gamma_\tau$. Thus, the matrix $\mtx{C}_{\tau}$ is used to convert, or translate, a set of MFS source strengths into a set of source strengths on the collocation points of $\Gamma_\tau$.

We now want to compress $\vct{q}_{\tau}$. Evaluating \cref{eq:out_skel} at the proxy points $\pzz_j$ on $\Psi_\tau$ and requiring that the equivalent sources $\hat{\vct{q}}_{\tau}$ reproduce the outgoing field there to precision $\varepsilon$, we get
\begin{equation}
    \mtx{B}_\tau(\colon,\hat{I}_{\tau}) \hat{\vct{q}}_{\tau} \stackrel{\varepsilon}{=} [u_\tau(\pzz_j)]_{j=1}^p = \mtx{B}_\tau \vct{q}_{\tau} \stackrel{\varepsilon}{=} \mtx{B}_\tau(\colon,\hat{I}_{\tau})\mtx{Z}_{\tau}^{*} \vct{q}_{\tau},
\end{equation}
where we used \cref{eq:B_fac} in the last step. Thus, we can set $\hat{\vct{q}}_{\tau} \stackrel{\varepsilon}{=} \mtx{Z}_{\tau}^{*} \vct{q}_{\tau}$. In other words, the $\mtx{Z}_\tau^*$ matrix can be used to compress the source strengths on $\Gamma_\tau$.

\subsection{Building the scattering operator}
\label{sec:build-scat}

We can now put everything together to find a relation between the compressed outgoing and incoming representations $\hat{\vct{q}}_{\tau}$ and $\hat{\vct{w}}_{\tau}$, respectively. The result is
\begin{equation}
\hat{\vct{q}}_{\tau} \stackrel{\varepsilon}{=}
\mtx{Z}_{\tau}^{*} \vct{q}_{\tau} =
\mtx{Z}_{\tau}^{*}\mtx{C}_\tau\vct{q}_{\tau}^{\text{mfs}} \approx
\mtx{Z}_{\tau}^{*}\mtx{C}_\tau\mtx{A}_\tau^{\dagger}\mtx{U}_{\tau}\hat{\vct{w}}_{\tau} =
\mtx{S}_{\tau}\hat{\vct{w}}_{\tau},
\end{equation}
where the local scattering matrix is defined as
\begin{equation}
\label{eq:scattering-matrix}
\mtx{S}_{\tau} = \mtx{Z}_{\tau}^{*}\mtx{C}_\tau \mtx{A}_\tau^{\dagger}\mtx{U}_{\tau}.
\end{equation}
Note that this is of the same form as \cref{eq:localnew} with $\mtx{V}_\tau^* = \mtx{Z}_\tau^* \mtx{C}_\tau$. 

Finally, let us note that inserting \cref{eq:scattering-matrix} into \cref{eq:mbody2} leads to
\begin{equation}
    \mtx{S}_{\tau}\hat{\vct{v}}_{\tau} = \mtx{Z}_{\tau}^{*}\mtx{C}_\tau \mtx{A}_\tau^{\dagger}\mtx{U}_{\tau}\hat{\vct{v}}_{\tau} = \mtx{Z}_{\tau}^{*}\mtx{C}_\tau \mtx{A}_\tau^{\dagger} \vct{v}_{\tau},
\end{equation}
which means that the effective charges $\mtx{S}_{\tau}\hat{\vct{v}}_{\tau}$ are computed per surface by applying $\mtx{Z}_{\tau}^{*}\mtx{C}_\tau \mtx{A}_\tau^{\dagger}$ to the global incoming field $\vct{v}_\tau$.

An alternative would be to select a subset of the MFS sources as skeleton points for the outgoing field, instead of translating the representation to the collocation skeletons. In this case, we would need to conduct skeletonizations on both collocation points and MFS sources. Then the local scattering matrix would follow the form of \cref{eq:localnew}, with independent matrices $\mtx{V}_\tau$ and $\mtx{U}_\tau$. 
This alternative method works, but we have found empirically that using the same points leads to smaller and more stable condition numbers as the resolution varies. 

\subsection{Reconstructing the full solution}
\label{sec:reconstruction}
The solved density $\hat{\vct{q}}$ is defined on the skeleton points on $\Gamma$. We can reconstruct the full solution either on $\Gamma_{\tau}$ or $\Gamma_{\tau}^{\text{mfs}}$, depending on where we want to evaluate the outgoing (scattered) field.

If reconstructed on the physical surface $\Gamma_\tau$, the full solution is
\begin{equation}
    \vct{q}_{\tau} = \mtx{C}_\tau \mtx{A}_\tau^{\dagger} \bigl(\vct{v}_{\tau} - \mtx{U}_{\tau} \hat{\mtx{G}}_{\tau,\tau^c} \hat{\vct{q}}_{\tau^c}\bigr),
\end{equation}
where $\hat{\mtx{G}}_{\tau,\tau^c}$ is the matrix whose entries
are the fundamental solution evaluated between skeleton points on
$\Gamma_\tau$ and $\Gamma \setminus \Gamma_\tau$, and $\hat{\vct{q}}_{\tau^c}$ is the solved density on the skeleton points of $\Gamma \setminus \Gamma_\tau$.

If reconstructed on the MFS surface $\Gamma_\tau^\text{mfs}$, the full solution is 
\begin{equation}
    \vct{q}_{\tau}^{\text{mfs}} = \mtx{A}_\tau^{\dagger} \bigl(\vct{v}_{\tau} - \mtx{U}_{\tau} \hat{\mtx{G}}_{\tau,\tau^c} \hat{\vct{q}}_{\tau^c}\bigr).
\end{equation}
An advantage of reconstructing on the MFS surface instead of the physical surface is that there is no singularity at $\Gamma_\tau$. The field can thus be evaluated at any point on the physical surface, and there is no need for special quadrature. Furthermore, the field can be evaluated in $\Omega$ arbitrarily close to $\Gamma_\tau$, which would also require special quadrature if a boundary integral equation were to be used.

Fast algorithms, like the fast multipole method, can also be applied in the reconstruction step, due to the nature of $\hat{\mtx{G}}$. The dominant cost in the reconstruction is then typically the application of $\mtx{A}_\tau^\dagger$ (which is completely local for each scatterer).

\section{Numerical experiments}
\label{sec:NumericalResults}

This section deals with some numerical experiments applying our solver to the Helmholtz problem, \textit{i.e.}, \cref{eq:BVP}.
We present two sets of numerical experiments: In \cref{sec:2d-results}, the focus is on 2D problems. We explore the convergence of our method for a number of small-scale 2D problems with differently shaped scatterers, including objects with cavities or corners. Furthermore, we demonstrate the performance of our solver in a 2D scaling problem where we increase the number of scatterers. Then, in \cref{sec:3d-results}, we consider problems in 3D, again starting with a few small-scale problems and concluding with a 3D scaling problem.

The numerical computations in this section were carried out using MATLAB. All timings are reported on a workstation with an Intel(R) Xeon(R) Gold 6248R processor (48 cores with 2 threads per core at 3.00\,GHz base frequency) and 1.5\,TB of RAM (the code will also run without issues on, \textit{e.g.}, a standard laptop). However, we point out that our code is not yet fully optimized, and that code performance is not the focus of this paper (for example, we do not report parallel performance in the form of, \textit{e.g.}, strong or weak scaling tests).

\subsection{Experiments in two dimensions}
\label{sec:2d-results}

\begin{example}
    We consider the exterior Helmholtz problem on a 2D ``starfish" geometry shown in \cref{fig: Starfish2D}. Here, each $\Gamma_\tau$ is a starfish-shaped boundary. The starfish shape is a smooth curve, parametrized by $(x_1,x_2) = (f(t) \cos t, f(t) \sin t)$, $f(t) = \frac{81}{101} - \frac{20}{101} \cos(5t)$, $t \in [0, 2\pi)$, that we discretize by points placed equidistantly in parameter space, as mentioned in \cref{sec:mfs-smooth-2d}.

    \begin{figure}[ht]
     \centering
     \includegraphics[width=\textwidth,trim=2cm 7cm 2cm 7cm,clip]{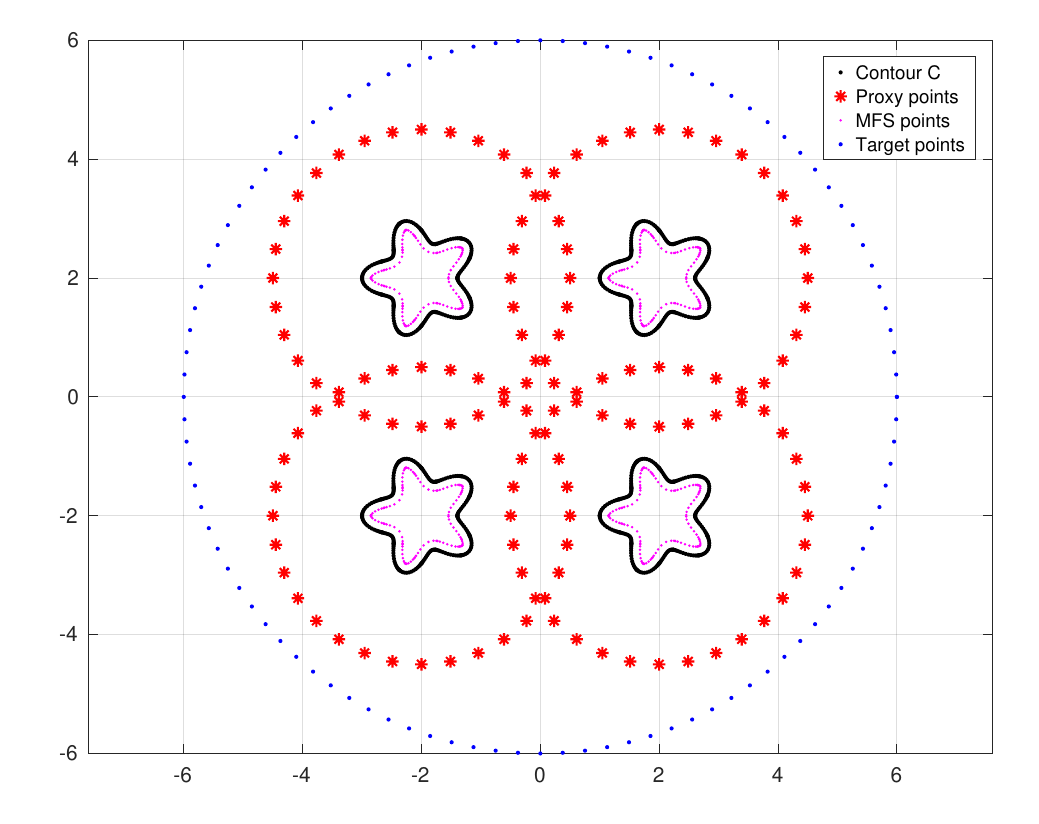}   \caption{The 2D 4-starfish geometry}
        \label{fig: Starfish2D}
    \end{figure} 

    We fix the number of scatterers to $T=4$ and vary the number of collocation points per scatterer. The externally applied incoming field is taken to be the plane wave $e^{i\kappa x_1}$. We run tests with wavenumbers $\kappa = 1, \pi, 10, 25$, so that each scatterer $\Gamma_{\tau}$ is approximately $0.3, 1, 3, 8$ wavelengths across (the diameter of a starfish is approximately 2; the wavelength is $\lambda = 2\pi/\kappa$). Results are presented in \cref{tab:Starfish2D-far-field,tab:Starfish2D-contours_pi,tab:Starfish2D-contours_10,tab:Starfish2D-contours_25}. We note the following quantities:
    \begin{center}
    \begin{tabular}{@{} l  l @{}}
        $N$  & The number of collocation points per scatterer (before skeletonization). \\
        $d$    & The MFS distance (distance between $\Gamma_\tau$ and $\Gamma_\tau^\text{mfs}$). \\
        $N_{\rm skel}$ & The number of skeleton points per scatterer. \\
        $K$ & The condition number of the global coefficient matrix $\mtx{I} + \mtx{S}\hat{\mtx{G}}$, \textit{cf.} \cref{eq:mbody2}. \\
        $E_{\rm far}$ & The relative error in the far-field potential (evaluated at blue dots in \cref{fig: Starfish2D}). \\
        $E_{\rm inc}$ & The relative error in the incoming field $w_\tau$. \\
        $\varepsilon$ & The skeletonization error tolerance.\\
        $m$ & The number of panels per scatterer.\\
        $p$ & The number of collocation points on each panel.\\
        $N_{\rm refine}$ & The number of dyadic refinements of panels around geometrical singularities.\\
    \end{tabular}
    \end{center}
    (The last three quantities are only applicable to non-smooth geometries and will appear in later examples.)

    $E_{\rm far}$ is the maximum relative error measured at the far-field target points shown as blue dots in Figure~\ref{fig: Starfish2D}. $E_{\rm inc}$ is the maximum relative error measured in the incoming field $w_\tau$ evaluated at fixed points placed uniformly on the physical contours $\Gamma_\tau$. Recall that the incoming field is the sum of the externally applied field and the scattered fields from the other scatterers, as defined in \cref{sec:introS} (in other words, the contribution from $\Gamma_\tau$ itself is excluded). The number of MFS sources is always half the number of collocation points, and the number of proxy points is selected as $p=N+1$. The reference solution in this example is computed using our method with $N=704$ collocation points per scatterer and $d = 0.08$.
    
    \begin{table}[ht]
        \centering\small
        \caption{2D 4-starfish, our method ($\kappa = 1$, $\varepsilon$ = $10^{-10}$)}
        \label{tab:Starfish2D-far-field}
        \begin{tabular}{cccccc}
        \toprule
        $N$ &
        $d$ &
        $N_{\rm skel}$ &
        $K$ &
        $E_{\rm far}$ &
        $E_{\rm inc}$ \\
        \midrule
        192 & 0.1 &  39   & 61.00  & $4.21\times 10^{-9}$ & $1.37\times 10^{-8}$ \\
        256 & 0.1 &  39   & 60.33  & $3.59\times 10^{-11}$ & $1.04\times 10^{-10}$ \\
        352 & 0.08 &  40   & 66.98  & $3.42\times 10^{-13}$ & $9.11\times 10^{-13}$ \\
        \bottomrule
        \end{tabular}
    \end{table}

    \begin{table}[ht]
        \centering\small
        \caption{2D 4-starfish, our method ($\kappa = \pi$, $\varepsilon$ = $10^{-10}$)}
        \label{tab:Starfish2D-contours_pi}
        \begin{tabular}{cccccc}
        \toprule
        $N$ &
        $d$ &
        $N_{\rm skel}$ &
        $K$ &
        $E_{\rm far}$ &
        $E_{\rm inc}$ \\
        \midrule
        192 & 0.1 &  42   & 9.05  & $9.01\times 10^{-9}$ & $9.46\times 10^{-9}$ \\
        256 & 0.1 &  42   & 10.33  & $6.10\times 10^{-11}$ & $6.65\times 10^{-11}$\\
        352 & 0.08 &  42   & 8.95  & $5.35\times 10^{-13}$ & $3.25\times 10^{-13}$\\
        \bottomrule
        \end{tabular}
    \end{table}

    \begin{table}[ht]
        \centering\small
        \caption{2D 4-starfish, our method ($\kappa = 10$, $\varepsilon$ = $10^{-10}$)}
        \label{tab:Starfish2D-contours_10}
        \begin{tabular}{cccccc}
        \toprule
        $N$ &
        $d$ &
        $N_{\rm skel}$ &
        $K$ &
        $E_{\rm far}$ &
        $E_{\rm inc}$ \\
        \midrule
        192 & 0.1 &  54   & 12.66  & $7.19\times 10^{-7}$ & $4.06\times 10^{-7}$ \\
        256 & 0.1 &  54   & 18.01  & $8.16\times 10^{-9}$ & $4.97\times 10^{-9}$\\
        352 & 0.08 &  54   & 27.73  & $2.45\times 10^{-11}$ & $2.45\times 10^{-11}$\\
        \bottomrule
        \end{tabular}
    \end{table}

    \begin{table}[ht]
        \centering\small
        \caption{2D 4-starfish, our method ($\kappa = 25$, $\varepsilon$ = $10^{-8}$)}
        \label{tab:Starfish2D-contours_25}
        \begin{tabular}{cccccc}
        \toprule
        $N$ &
        $d$ &
        $N_{\rm skel}$ &
        $K$ &
        $E_{\rm far}$ &
        $E_{\rm inc}$ \\
        \midrule
        192 & 0.1 &  85   & 41.88 & $7.30\times 10^{-4}$  & $6.91\times 10^{-4}$ \\
        256 & 0.1 &  85   & 31.24 & $6.05\times 10^{-6}$ & $6.43\times 10^{-6}$\\
        352 & 0.08 &  85   & 29.27 & $5.31\times 10^{-9}$ & $5.70\times 10^{-9}$\\
        \bottomrule
        \end{tabular}
    \end{table}

    \begin{table}[ht]
        \centering\small
        \caption{2D 4-starfish, BIE method ($\kappa = 25, p = 16$)}
        \label{tab:Starfish2D_25_BIE}
        \begin{tabular}{ccccc}
        \toprule
        $N$ &
        $m$ &
        $K$ &
        $E_{\rm far}$ &
        $E_{\rm inc}$ \\
        \midrule
        \num{192} & \num{12} & \num{20.57} & \num{5.43e-06} & \num{3.93e-06} \\
        \num{256} & \num{16} & \num{20.60} & \num{4.66e-08} & \num{1.19e-08} \\
        \num{352} & \num{22} & \num{20.26} & \num{1.31e-09} & \num{8.34e-10} \\
        \bottomrule
        \end{tabular}
    \end{table}

    \Cref{tab:Starfish2D-far-field,tab:Starfish2D-contours_pi,tab:Starfish2D-contours_10,tab:Starfish2D-contours_25} show that $E_{\rm far}$ is basically at the same level as $E_{\rm inc}$, that they both decay spectrally when the resolution increases, and that they increase moderately with the wavenumber (for fixed resolution). In the remaining examples, we will only use $E_{\rm inc}$ to report the accuracy of the solver, and perform experiments only with the highest wavenumber considered here ($\kappa=25$).

    Table \ref{tab:Starfish2D_25_BIE} shows the results of a BIE solver using the \textit{chunkIE} package \citep{chunkie}, with similar resolution. The BIE solver uses a discretization based on Gauss--Legendre panels, with $m$ panels per scatterer and $p$ points per panel. Comparison between Tables \ref{tab:Starfish2D-contours_25} and \ref{tab:Starfish2D_25_BIE} suggests that the performance of our stable MFS solver is comparable with the BIE solver for these smooth starfish geometries.
\end{example}

\begin{example}
    We now consider a C-shaped scatterer with a cavity, shown in \cref{fig: SmoothLune2D}. The C-shape consists of two circular arcs (one inner and one outer) joined by two semicircular caps. Each segment is discretized using $m_{\rm seg}$ panels. Dyadic refinement is needed at the joints of these segments, since the second derivative of the geometry is discontinuous there. We repeat the experiment of the previous example for this geometry with $\kappa = 25$, so that $\Gamma_{\tau}$ is approximately 8 wavelengths across. The convergence results are presented in \cref{tab:SmoothLune2D_25_2}. Here, $N_0$ is the base resolution (before dyadic refinement), and $N$ is the resolution after refinement. The reference solution is computed by the same method with parameters $m_{\rm seg} = 256, p = 16, d = 0.05, N_{\rm refine} = 10, \varepsilon = 10^{-10}$.

    \begin{figure}[ht]
     \centering
     \includegraphics[width=\textwidth,trim=2cm 7cm 2cm 7cm,clip]{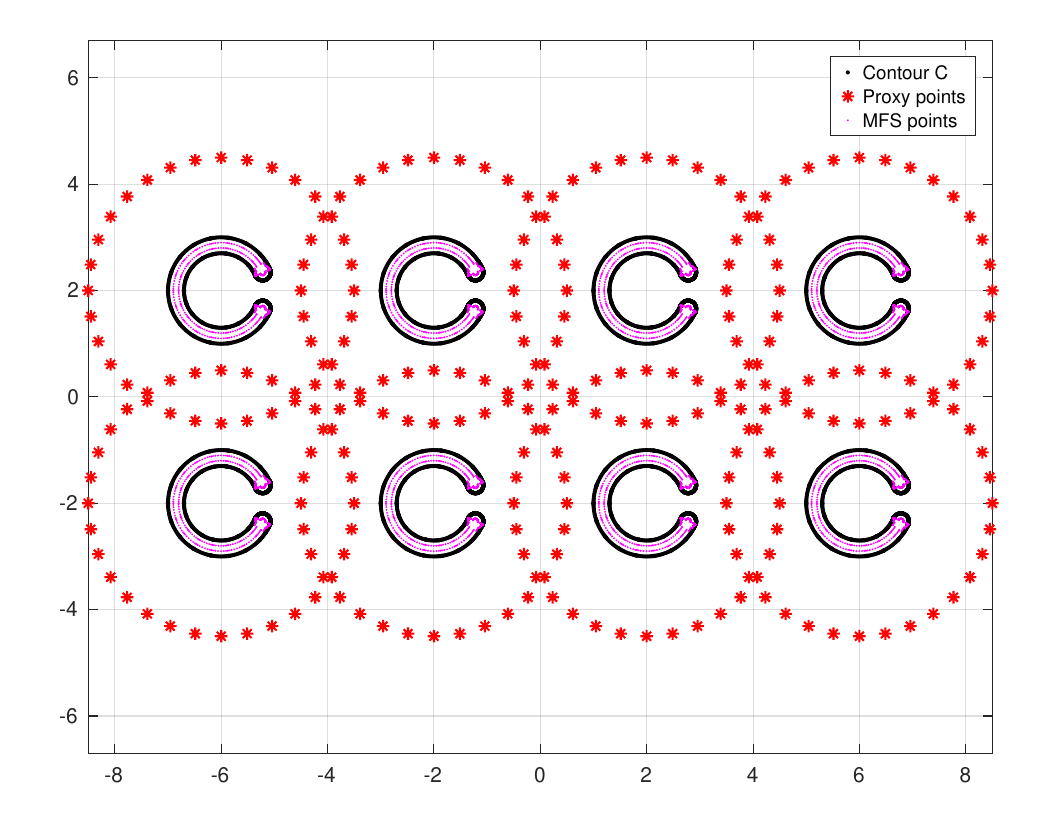}
        \caption{The 2D 8-cavity geometry}
        \label{fig: SmoothLune2D}
    \end{figure} 

    \begin{table}[ht]
        \centering\small
        \caption{2D 8-cavity, our method ($\kappa = 25, \varepsilon = 10^{-10}$, $p=16$), with externally applied incoming field $e^{i\kappa x_1}$}
        \label{tab:SmoothLune2D_25_2}
        \begin{tabular}{ccccccccc}
        \toprule
        $N_0$ &
        $m_{\rm seg}$ &
        $m$ &
        $N_{\rm refine}$ &
        $N$ &
        $d$ &
        $N_{\rm skel}$ &
        $K$ &
        $E_{\rm inc}$ \\
        \midrule
        512 & \num{8} & \num{32} & \num{5} & 1152 & \num{0.10} & \num{97} & \num{114.72} & \num{6.97e-03} \\
        1024 & \num{16} & \num{64} & \num{5} & 1664 & \num{0.10} & \num{98} & \num{335.07} & \num{2.08e-06} \\
        2048 & \num{32} & \num{128} & \num{5} & 2688 & \num{0.10} & \num{98} & \num{215.79} & \num{3.84e-10} \\
        4096 & \num{64} & \num{256} & \num{10} & 5376 & \num{0.10} & \num{98} & \num{264.55} & \num{8.82e-11} \\
        8192 & \num{128} & \num{512} & \num{10} & 9472 & \num{0.05} & \num{98} & \num{176.63} & \num{6.78e-11} \\
        \bottomrule
        \end{tabular}
    \end{table}

\begin{table}[ht]
        \centering\small
        \caption{2D 8-cavity, BIE method ($\kappa = 25, p = 16$)}
        \label{tab:SmoothLune2D_25_BIE}
        \begin{tabular}{ccccc}
        \toprule
        $N_0$ &
        $m_{\rm seg}$ &
        $m$ &
        $K$ &
        $E_{\rm inc}$ \\
        \midrule
        \num{704} & \num{8} & \num{44} & \num{333.07} & \num{1.90e-07} \\
        \num{1216} & \num{16} & \num{76} & \num{324.98} & \num{6.91e-11} \\
        \num{2240} & \num{32} & \num{140} & \num{315.88} & \num{1.19e-10} \\
        \num{4288} & \num{64} & \num{268} & \num{310.49} & \num{2.14e-10} \\
        \num{8384} & \num{128} & \num{524} & \num{306.71} & \num{3.94e-10} \\
        \bottomrule
        \end{tabular}
    \end{table}
    
    The results show that our stable MFS solver behaves well even when the geometry has cavities. Again, the results are comparable with the BIE solver, shown in \cref{tab:SmoothLune2D_25_BIE}. The BIE solver \textit{chunkIE} selects $m_{\rm seg}$ panels on the small semicircular caps, and $m_{\rm seg}+6$ panels on the inner and outer circular arcs.
    It should be noted that, in our solver, $N$ is the number of unknowns only in the local problem (when constructing the scattering matrix or reconstructing the full solution). In the global system, the number of unknowns is $N_\mathrm{skel}$ times the number of scatterers, which is significantly smaller. For the BIE solver, the number of unknowns is $N_0$ times the number of scatterers.
\end{example}

\begin{example}
    Here, we turn to a geometry with corners, specifically teardrop-shaped structures. The configuration consists of 8 scatterers (see \cref{fig: Teardrop2D}). Each teardrop is parametrized as $(x_1,x_2) = (\frac{2}{\pi^2} t^2 - \frac{4}{\pi} t + 1, \frac{2}{\pi^3} t^3 - \frac{6}{\pi^2} t^2 + \frac{4}{\pi} t)$, $t\in [0, 2\pi)$. We discretize each scatterer with $m$ panels and $p=16$ Gauss--Legendre nodes on each panel. To deal with the singularity at the corner ($t=0$), the panels next to the corner are dyadically refined, as mentioned in \cref{sec:mfs-nonsmooth-2d}.
    
    We set the wavenumber $\kappa = 25$ and skeletonization tolerance $\varepsilon = 10^{-10}$. In Table \ref{tab:Teardrop2D_25}, we fix $N_{\mathrm{refine}} = 20$ dyadically refined panels in the vicinity of corners, and increase the number of base panels $m$ per scatterer. The reference solution is computed by the same method with parameters $m = 128, p = 16, d = 0.1, N_{\rm refine} = 50, \varepsilon = 10^{-10}$.

   \begin{figure}[ht]
     \centering
     \includegraphics[width=\textwidth,trim=2cm 7cm 2cm 7cm,clip]{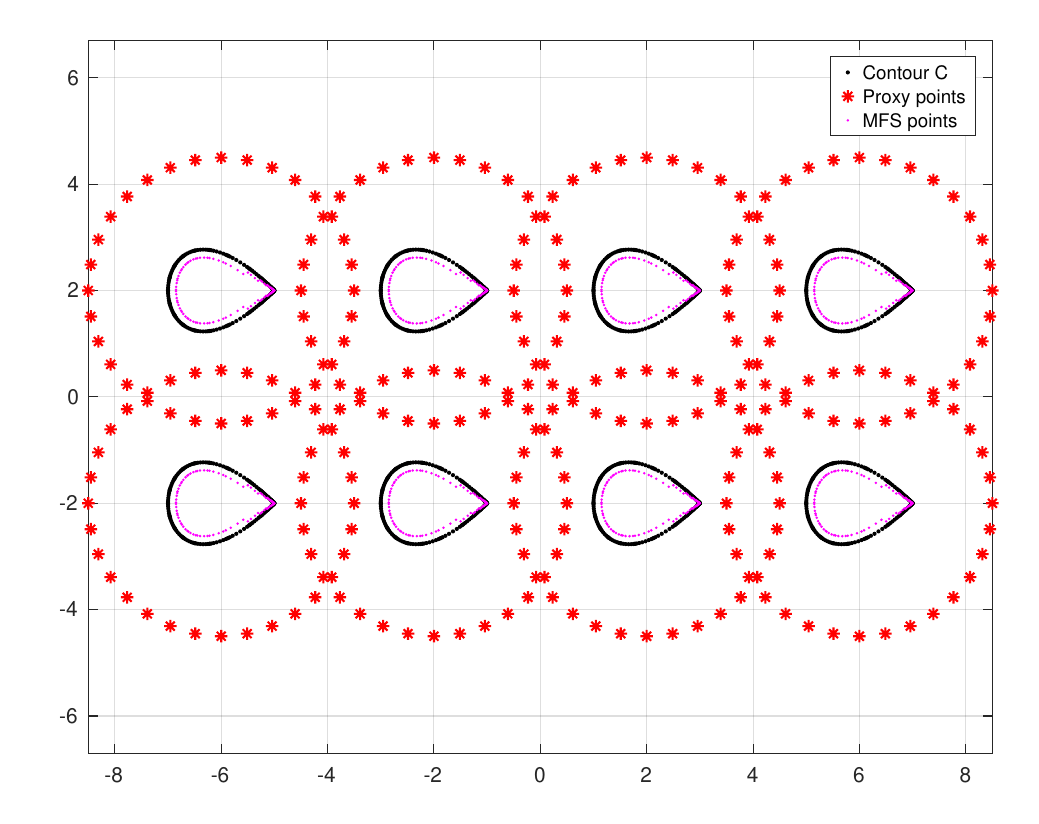}
        \caption{The 2D 8-teardrop geometry}
        \label{fig: Teardrop2D}
    \end{figure} 

    \begin{table}[ht]
        \centering\small
        \caption{2D 8-teardrop, our method ($\kappa = 25, \varepsilon = 10^{-10}$, $p=16$, $N_{\rm refine} = 20$)}
        \label{tab:Teardrop2D_25}
        \begin{tabular}{ccccccc}
        \toprule
        $N_0$ &
        $m$ &
        $N$ &
        $d$ &
        $N_{\rm skel}$ &
        $K$ &
        $E_{\rm inc}$ \\
        \midrule
        64 & 4 & 704 & 0.25 & 64   & 51.26  & $4.16\times 10^{-2}$ \\
        128 & 8 & 768 & 0.25 & 88   & 54.10  & $1.17\times 10^{-4}$ \\
        256 & 16 & 896 & 0.25 & 88   & 41.39  & $3.75\times 10^{-8}$ \\
        512 & 32 & 1152 & 0.2 & 89  & 41.58  & $1.43\times 10^{-9}$ \\
        1024 & 64 & 1664 & 0.1 & 88   & 38.09  & $2.49\times 10^{-10}$ \\
        2048 & 128 & 2688 & 0.1 & 88   & 39.16  & $2.73\times 10^{-10}$ \\
        \bottomrule
        \end{tabular}
    \end{table}

    \begin{table}[ht]
        \centering\small
        \caption{2D 8-teardrops, BIE method ($\kappa = 25, p = 16$)}
        \label{tab:Teardrop2D_25_BIE}
        \begin{tabular}{ccccc}
        \toprule
        $N_0$ &
        $m$ &
        $K$ &
        $E_{\rm inc}$ \\
        \midrule
        \num{64} & \num{4} & \num{626.42} & \num[retain-zero-exponent = true]{3.68e+00} \\
        \num{128} & \num{8} & \num{40.95} & \num{1.16e-04} \\
        \num{256} & \num{16} & \num{40.47} & \num{1.35e-05} \\
        \num{512} & \num{32} & \num{40.27} & \num{5.06e-06} \\
        \num{1024} & \num{64} & \num{39.78} & \num{1.95e-06} \\
        \num{2048} & \num{128} & \num{39.67} & \num{7.64e-07} \\
        \num{4096} & \num{256} & \num{39.68} & \num{3.01e-07} \\
        \bottomrule
        \end{tabular}
    \end{table}

    Table \ref{tab:Teardrop2D_25_BIE} shows the results of a BIE solver, with similar setting of discretization. Comparison between Tables \ref{tab:Teardrop2D_25} and \ref{tab:Teardrop2D_25_BIE} suggests that the performance of our stable MFS solver is comparable with the BIE solver also for the teardrop-shaped geometry, if not slightly better. The BIE solver (\textit{chunkIE}) also uses dyadic refinement around the corner, in the form of recursively compressed inverse preconditioning (RCIP), which does not affect the number of unknowns. Still, the convergence of the BIE solver is somewhat slow in this case. However, this may be due to implementation details and does not imply that BIE-based methods are generally inferior to our method for domains with corners.
\end{example}

\begin{example}
\label{ex:2d-results-scaling}

As a final two-dimensional experiment, we consider a scaling problem with an increasing number of scatterers. The geometry consists of starfishes, C-shaped cavities, teardrops and rods (our rods are thin rectangles with semicircular caps attached). The number of scatterers is scaled up from $T=4$ to $T=2048$; an example of the geometry is shown in \cref{fig:visualization2d}, where the solution (incoming + scattered field) is also visualized. The wavenumber is $\kappa=25$ (each scatterer is about 8 wavelengths across; it should be pointed out that this is a rather high wavenumber compared to related work such as, \textit{e.g.}, \cite{2022_stein_barnett} which was mentioned in \cref{sec:introduction}) and the external incoming field is $e^{i \kappa x_1}$. The global linear system, \textit{i.e.}, \cref{eq:mbody2}, is solved iteratively using GMRES (MATLAB's built-in \texttt{gmres} function), and the matrix-vector products (matvecs) are accelerated by the FMM (the \texttt{fmm2d} package \citep{fmm2d}). 

\begin{figure}[htp]
  \centering
  \includegraphics[width=0.33\textwidth]{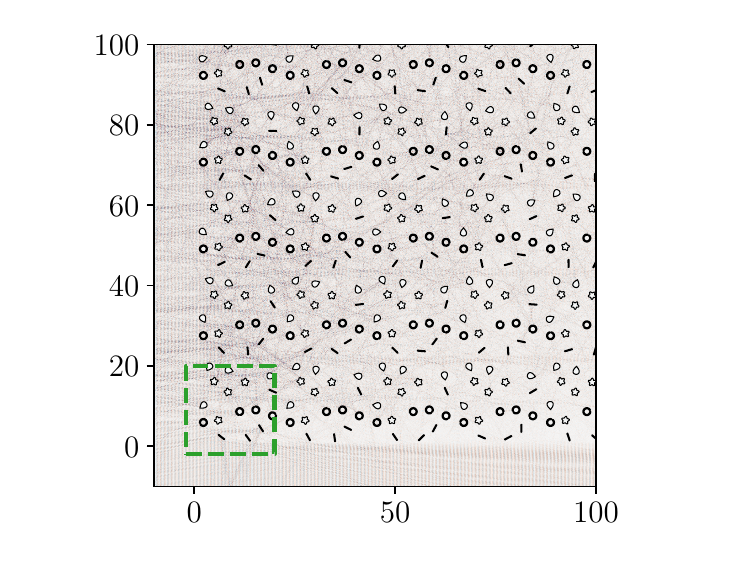}%
  \includegraphics[width=0.33\textwidth]{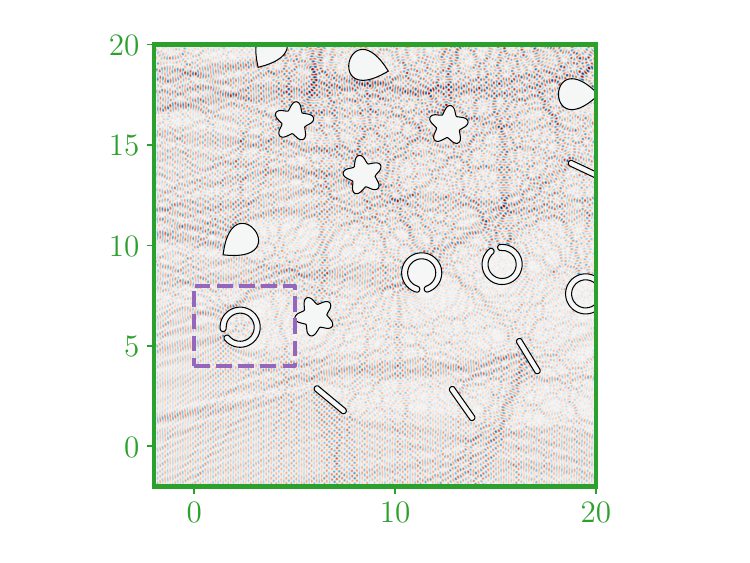}%
  \includegraphics[width=0.33\textwidth]{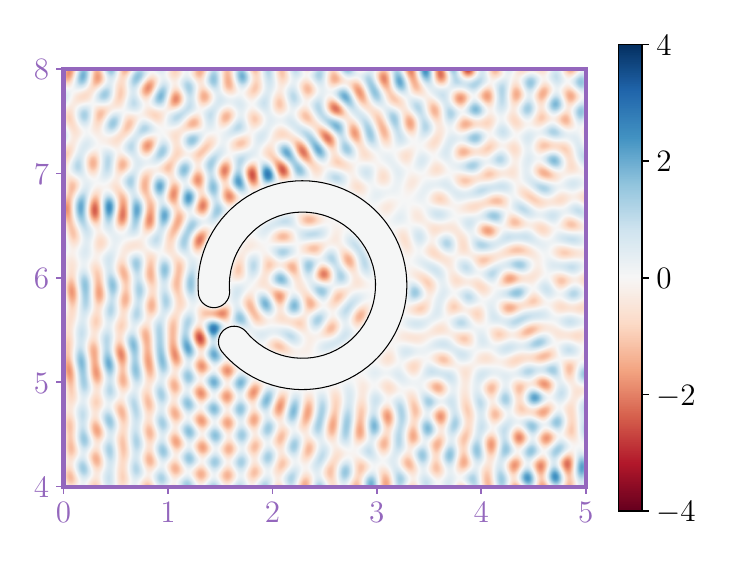}%
  \caption{Imaginary part of total field (incoming + scattered) for a problem involving $T=1024$ scatterers ($\kappa=25$, incoming field $e^{i \kappa x_1}$), at three different zoom levels}
  \label{fig:visualization2d}
\end{figure}

Results for parameter choices corresponding to two different error levels%
\footnote{The number of collocation points per scatterer ($N$), including points from dyadic refinement, is shown for each error level in the table (starfishes/C-shapes/teardrops/rods). Other parameters are, for error level $10^{-5}$: MFS distance ($d$) 0.08/0.12/0.30/$\frac{2}{15}$, number of dyadic refinements ($N_{\rm refine}$) 0/0/10/1. For error level $10^{-9}$: MFS distance ($d$) 0.08/0.10/0.15/0.073, number of dyadic refinements ($N_{\rm refine}$) 0/4/25/4. The number of points per Gauss--Legendre panel is always $p=16$ for all non-smooth shapes (C-shapes, teardrops, rods). For rods, approximately $\frac{4}{11}$ of the panels are placed on each straight line segment, and $\frac{3}{22}$ on each semicircular cap. For C-shapes, approximately $\frac{2}{5}$ of the panels are placed on the outer circular arc, $\frac{3}{10}$ on the inner circular arc, and $\frac{3}{20}$ on each semicircular cap. These ratios refer to base panels ($m$), before dyadic refinement.}
are presented in \cref{tab:scaling_2d_final}. The error tolerances for skeletonization, GMRES and FMM are set to the stated error level divided by 10. The columns of \cref{tab:scaling_2d_final} show (from left to right): ($T$) the number of scatterers in the problem, ($n_\text{matvec}$) the number of matrix-vector products used by GMRES to solve the global linear system, ($N_\text{tot}$) the total number of collocation points on all scatterers (before skeletonization), ($N_\text{skel,tot}$) the total number of skeleton points on all scatterers; the runtime subdivided into ($t_\text{local}$) local computations (this includes computing the right-hand side of \cref{eq:mbody2} given the incoming field, and also reconstructing the full solution $\vct{q}$ after the equivalent sources $\hat{\vct{q}}$ have been determined, as described in \cref{sec:reconstruction}), ($t_\text{solve}$) total solve time for the iterative solver, ($t_\text{matvec}$) average time per matvec ($=t_\text{solve}/n_\text{matvec}$); and ($E_\text{inc}$) the maximum relative error of the incoming field.

The time spent on constructing the scattering matrices $\mtx{S}_\tau$, including skeletonization and computing the matrices $\mtx{A}_\tau^\dagger$ and $\mtx{C}_\tau$, is shown as ``precomp.\ time'' in the table. Since $\mtx{S}_\tau$ is identical for objects of the same shape (modulo a phase factor when rotated \citep{2015_lai_kobayashi_barnett}), only four scattering matrices need to be constructed (one for each shape), and this time is thus independent of the number of scatterers $T$.

\begin{table}[htp]
  \centering\small
  \caption{Scaling problem with a mix of 2D scatterers
  (starfishes, C-shapes, teardrops, rods) and wavenumber
  $\kappa=25$}
  \label{tab:scaling_2d_final}
  \vspace*{5pt}
    \sisetup{table-format=1.2e-1,table-auto-round=true,tight-spacing=true}
    \setlength{\tabcolsep}{3pt}
    \begin{tabular}{cccc|S[retain-zero-exponent=true]S[retain-zero-exponent=true]S[retain-zero-exponent=true]|S}
    \toprule
    \multicolumn{8}{c}{\textit{Error level $10^{-5}$ ($250/608/464/352$ points per scatterer, precomp.\ time \SI{0.33}{s})}} \\
    \midrule
    \multicolumn{4}{c|}{\bf Number of} & \multicolumn{3}{c|}{\bf Time [s]} & \multicolumn{1}{c}{\bf Error} \\
    scatterers & matvecs & points & skeletons & {local comp.} & {solve} & {per matvec} & {incoming} \\
    $T$ & $n_\text{matvec}$ & $N_\text{tot}$ & $N_\text{skel,tot}$ & {$t_\text{local}$} & {$t_\text{solve}$} & {$t_\text{matvec}$} & {$E_\text{inc}$} \\
    \midrule
    \num{4} & \num{19} & \num{1674} & \num{294} & 2.965700e-02 & 3.211800e-01 & 1.690421e-02 & 4.392240e-06 \\
    \num{8} & \num{31} & \num{3348} & \num{588} & 3.251200e-02 & 5.560160e-01 & 1.793600e-02 & 5.419470e-06 \\
    \num{16} & \num{49} & \num{6696} & \num{1176} & 3.859200e-02 & 1.108513e+00 & 2.262271e-02 & 1.187150e-05 \\
    \num{32} & \num{74} & \num{13392} & \num{2352} & 6.042500e-02 & 2.436980e+00 & 3.293216e-02 & 1.321000e-05 \\
    \num{64} & \num{134} & \num{26784} & \num{4704} & 1.042230e-01 & 5.856494e+00 & 4.370518e-02 & 2.261290e-05 \\
    \num{128} & \num{234} & \num{53568} & \num{9408} & 1.819270e-01 & 1.885466e+01 & 8.057545e-02 & 2.733270e-05 \\
    \num{256} & \num{524} & \num{107136} & \num{18816} & 2.874970e-01 & 7.052500e+01 & 1.345897e-01 & 2.868460e-05 \\
    \num{512} & \num{943} & \num{214272} & \num{37632} & 5.218240e-01 & 2.566480e+02 & 2.721612e-01 & 3.072300e-05 \\
    \num{1024} & \num{2377} & \num{428544} & \num{75264} & 1.120747e+00 & 1.701349e+03 & 7.157547e-01 & 8.060860e-05 \\
    \num{2048} & \num{4304} & \num{857088} & \num{150528} & 2.073660e+00 & 7.419695e+03 & 1.723907e+00 & 9.148540e-05 \\
    \bottomrule
    \toprule
    \multicolumn{8}{c}{\textit{Error level $10^{-9}$ ($404/2032/1504/1472$ points per scatterer, precomp.\ time \SI{2.77}{s})}} \\
    \midrule
    \multicolumn{4}{c|}{\bf Number of} & \multicolumn{3}{c|}{\bf Time [s]} & \multicolumn{1}{c}{\bf Error} \\
    scatterers & matvecs & points & skeletons & {local comp.} & {solve} & {per matvec} & {incoming} \\
    $T$ & $n_\text{matvec}$ & $N_\text{tot}$ & $N_\text{skel,tot}$ & {$t_\text{local}$} & {$t_\text{solve}$} & {$t_\text{matvec}$} & {$E_\text{inc}$} \\
    \midrule
    \num{4} & \num{25} & \num{5412} & \num{350} & 5.195800e-02 & 4.756640e-01 & 1.902656e-02 & 2.018680e-10 \\
    \num{8} & \num{46} & \num{10824} & \num{700} & 8.663100e-02 & 9.587830e-01 & 2.084311e-02 & 6.734420e-10 \\
    \num{16} & \num{78} & \num{21648} & \num{1404} & 1.481780e-01 & 1.950099e+00 & 2.500127e-02 & 3.094570e-10 \\
    \num{32} & \num{122} & \num{43296} & \num{2808} & 2.628550e-01 & 4.171372e+00 & 3.419157e-02 & 5.202660e-10 \\
    \num{64} & \num{221} & \num{86592} & \num{5616} & 5.339660e-01 & 1.150647e+01 & 5.206546e-02 & 1.473340e-09 \\
    \num{128} & \num{394} & \num{173184} & \num{11232} & 9.982730e-01 & 3.780167e+01 & 9.594332e-02 & 1.424450e-09 \\
    \num{256} & \num{892} & \num{346368} & \num{22464} & 1.984270e+00 & 1.525400e+02 & 1.710090e-01 & 1.809390e-09 \\
    \num{512} & \num{1638} & \num{692736} & \num{44928} & 3.946838e+00 & 6.305600e+02 & 3.849573e-01 & 2.162210e-09 \\
    \num{1024} & \num{4174} & \num{1385472} & \num{89856} & 7.854652e+00 & 4.607610e+03 & 1.103883e+00 & 4.943200e-09 \\
    \num{2048} & \num{7666} & \num{2770944} & \num{179712} & 1.631862e+01 & 2.300042e+04 & 3.000316e+00 & 7.713370e-09 \\
    \bottomrule
    \end{tabular}
\end{table}

We again point out that the global linear system is of size $N_\text{skel,tot} \times N_\text{skel,tot}$, and that $N_\text{tot}$ does not appear as the size of any system to be solved in this formulation. The reconstructed full density $\vct{q}$ is a vector of size $N_\text{tot}$, but it is computed by $T$ local solves of size $N \times N$. To quantify the compression accomplished by skeletonization, we can define the \emph{space saving} as $\eta = 1 - N_\text{skel,tot}/N_\text{tot}$, which is a number between 0 and 1 (higher is better). In \cref{tab:scaling_2d_final}, the space saving is $\eta=0.82$ for error level $10^{-5}$, and $\eta=0.94$ for error level $10^{-9}$.

The table shows that the time for local computations $t_\text{local}$ scales linearly in the number of scatterers $T$. The time per matvec $t_\text{matvec}$ is also linear in $T$ (due to the FMM), although a superlinear growth can be seen for $T \geq 1024$ in this test. This is due to increasing internal overhead in the \texttt{gmres} function itself as the number of iterations grows large. Potential ways to mitigate this include using a restarting strategy (which could reduce the cost per iteration but increase the total number of iterations), preconditioning (to reduce the number of iterations), or a more optimized implementation of GMRES.

The number of matvecs $n_\text{matvec}$ required by the iterative solver depends on both the number of scatterers $T$ and the GMRES error tolerance, but not on other parameters such as the number of collocation points or skeletonization tolerance, as shown in \cref{fig:gmres-tol-variation}. Each curve corresponds
  to a different GMRES tolerance ($10^{-10}$, $10^{-8}$, $10^{-6}$,
  $10^{-4}$); all other parameters are set
  according to the error level. Colored curves with circular
  markers are error level $10^{-9}$, black triangles are error
  level $10^{-5}$. The number of matvecs is independent of the
  error level, and depends only on the number of scatterers and
  the GMRES tolerance.
This reflects that the conditioning of the global linear system is independent of the discretization, similar to a well-conditioned BIE formulation.

The number of matvecs grows approximately linearly with $T$. This reflects the conditioning of the underlying partial differential equation, and is observed also in, \textit{e.g.}, a well-conditioned BIE formulation. As an example, we record in \cref{tab:scaling_2d_bie} the number of matvecs required to solve a similar Helmholtz problem involving (for simplicity) only starfish scatterers, with the same wavenumber, incoming field, and average concentration of scatterers, using a BIE solver. The number of matvecs can be compared with the lower half of \cref{tab:scaling_2d_final} (error level $10^{-9}$). The condition number is measured in the 2-norm. The BIE solver uses a combined field formulation with 16th-order Zeta correction \citep{2020_martinsson_zetacorrection} for the singular quadrature, and no compression or preconditioning. Again, both the number of matvecs and the condition number grow approximately linearly with the number of scatterers $T$.

\begin{figure}[htp]
  \centering
  \includegraphics[width=0.5\textwidth]{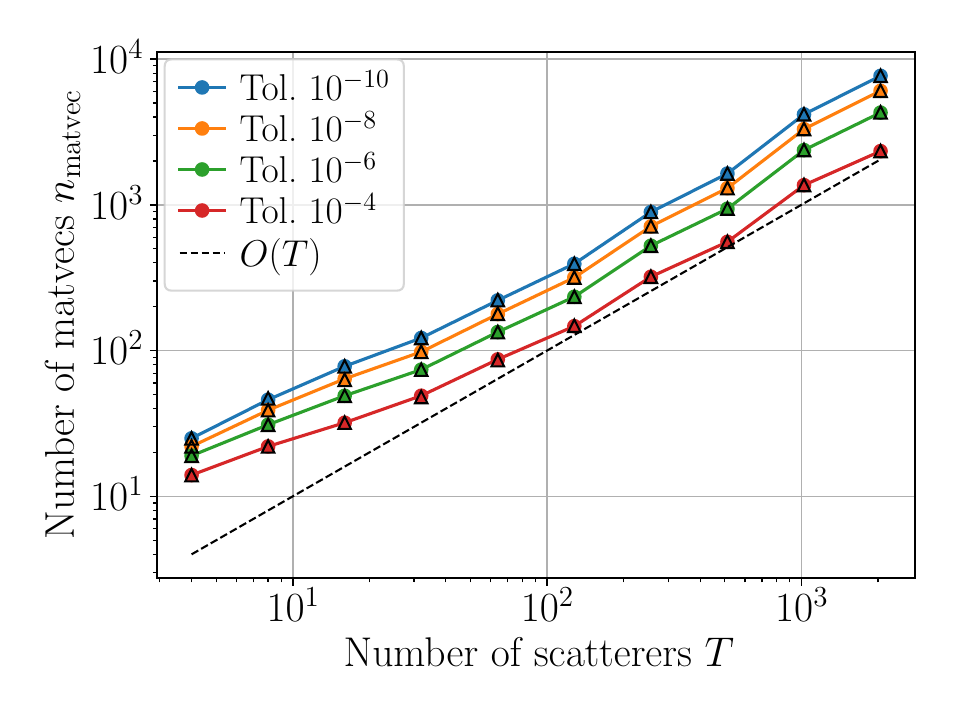}
  \caption{The number of matvecs $n_\text{matvec}$ during the solve stage as a function of the number of scatterers $T$ (one curve per GMRES tolerance)}
  \label{fig:gmres-tol-variation}
\end{figure}

\begin{table}[htp]
  \centering\small
  \caption{BIE scaling results for starfish scatterers, wavenumber $\kappa=25$, and GMRES error tolerance $10^{-10}$}
  \label{tab:scaling_2d_bie}
  \vspace*{5pt}
    \sisetup{table-format=3.1,table-auto-round=true,tight-spacing=true}
    \begin{tabular}{ccS}
    \toprule
    Number of scatterers & Number of matvecs & {Condition number} \\
    \midrule
    \num{4} & \num{47} & 15.2999 \\
    \num{8} & \num{69} & 23.5328 \\
    \num{16} & \num{133} & 65.3327 \\
    \num{32} & \num{201} & 111.554 \\
    \num{64} & \num{385} & 182.496 \\
    \bottomrule
    \end{tabular}
\end{table}

As can be seen in the rightmost column of \cref{tab:scaling_2d_final}, the observed error $E_\text{inc}$ grows slightly with the number of scatterers $T$. In order to keep the error fixed as $T$ increases, it is necessary to slightly increase the number of collocation points $N$, and also slightly reduce the GMRES tolerance. The other error tolerances (\textit{i.e.}, skeletonization and FMM) can be kept fixed. Note that $N$ mainly affects the accuracy of the local scattering matrix, as defined by \cref{eq:scattering-matrix}; it may be that the computation of this matrix requires higher resolution as an increase in $T$ leads to a more complicated incoming field from the point of view of an individual scatterer.

The singular values of the global coefficient matrix $\mtx{I} + \mtx{S} \hat{\mtx{G}}$ exhibit interesting behavior as the number of scatterers $T$ varies. For the Laplace problem (corresponding to $\kappa$ = 0), there is a clear gap between the first $T$ largest singular values and the rest. To be specific, the number of singular values larger than about $O(10^1)$ is exactly $T$; the $(T+1)\text{-th}$ singular value drops dramatically to $1+\varepsilon$ for some small $\varepsilon>0$. For the Helmholtz problem, there is no gap in the same sense, but the singular values smoothly decay to 1, with a rate depending on $T$. This is likely related to the number of matvecs $n_\text{matvec}$ depending on $T$, but a deeper analysis is outside the scope of this paper, and left for future work.

\end{example}

\subsection{Experiments in three dimensions}
\label{sec:3d-results}

\begin{example}
    We move on to three-dimensional experiments, and first consider geometries consisting of tori and ellipsoids separately. The configuration consists of 2 scatterers, as shown in \cref{fig: Torus3D,fig: Ellip3D}. The skeletonization tolerance $\varepsilon$ is fixed to $10^{-6}$, and we vary the number of collocation points $N$. The incoming field is a plane wave $e^{i\kappa x_1}$. The reference solution in this example is computed using the FMM3DBIE package \citep{fmm3dbie}. Results are shown in \cref{tab:EllipsoidTorus3D_5_10_3}, which suggests that the accuracy of our stable MFS solver is comparable with a BIE. $E_{\rm inc}$ decays spectrally with increasing resolution, and increases moderately with the wavenumber $\kappa$. Ellipsoids reach a smaller error than tori for the same $N$.
    
    \begin{figure}[ht]
     \centering
     \includegraphics[width=0.8\textwidth,trim=10cm 2cm 8cm 0.5cm,clip]{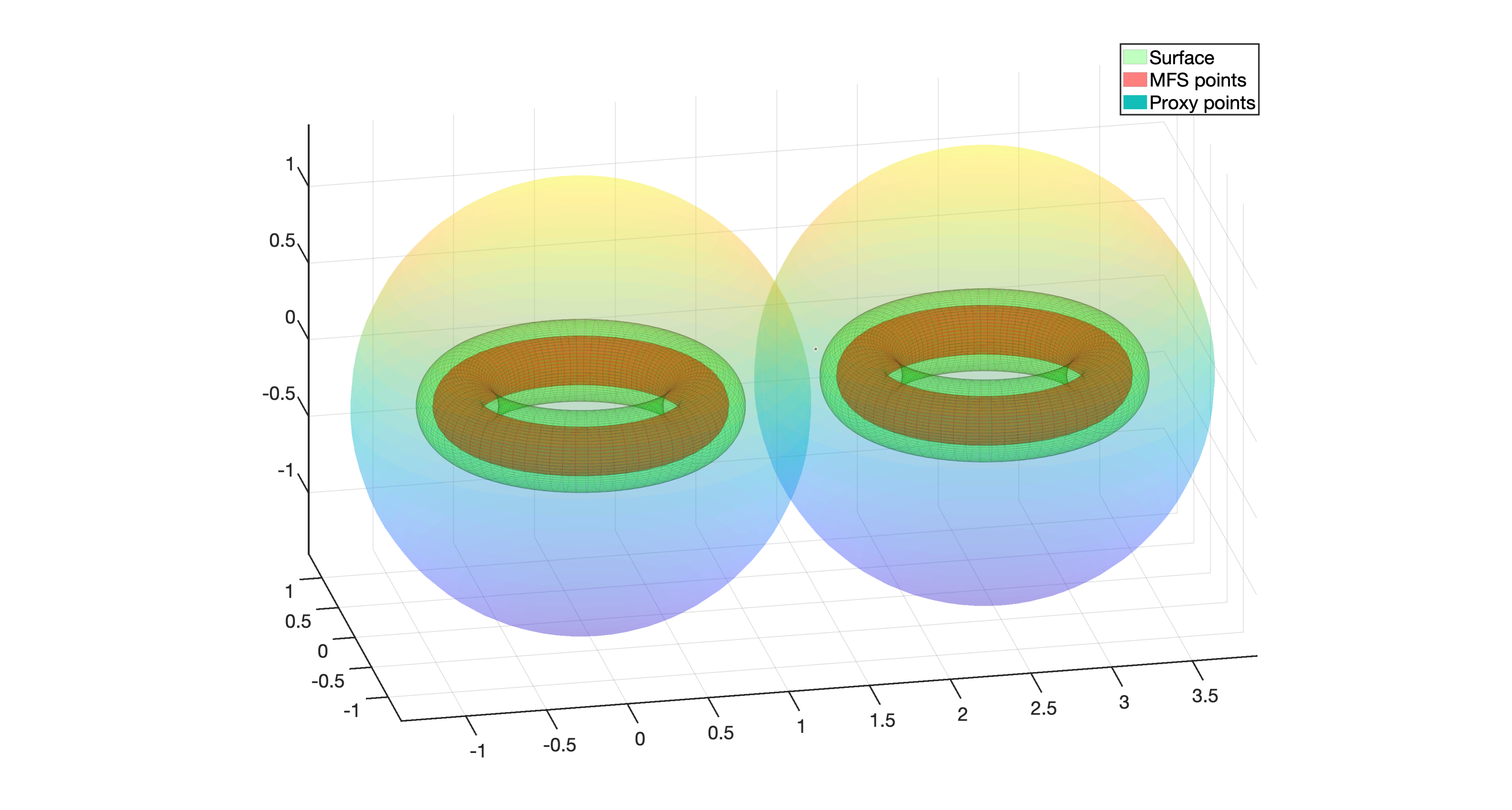}
        \caption{3D Torus geometry}
        \label{fig: Torus3D}
    \end{figure} 

    \begin{figure}[ht]
     \centering
     \includegraphics[width=0.8\textwidth,trim=10cm 2cm 8cm 0.5cm,clip]{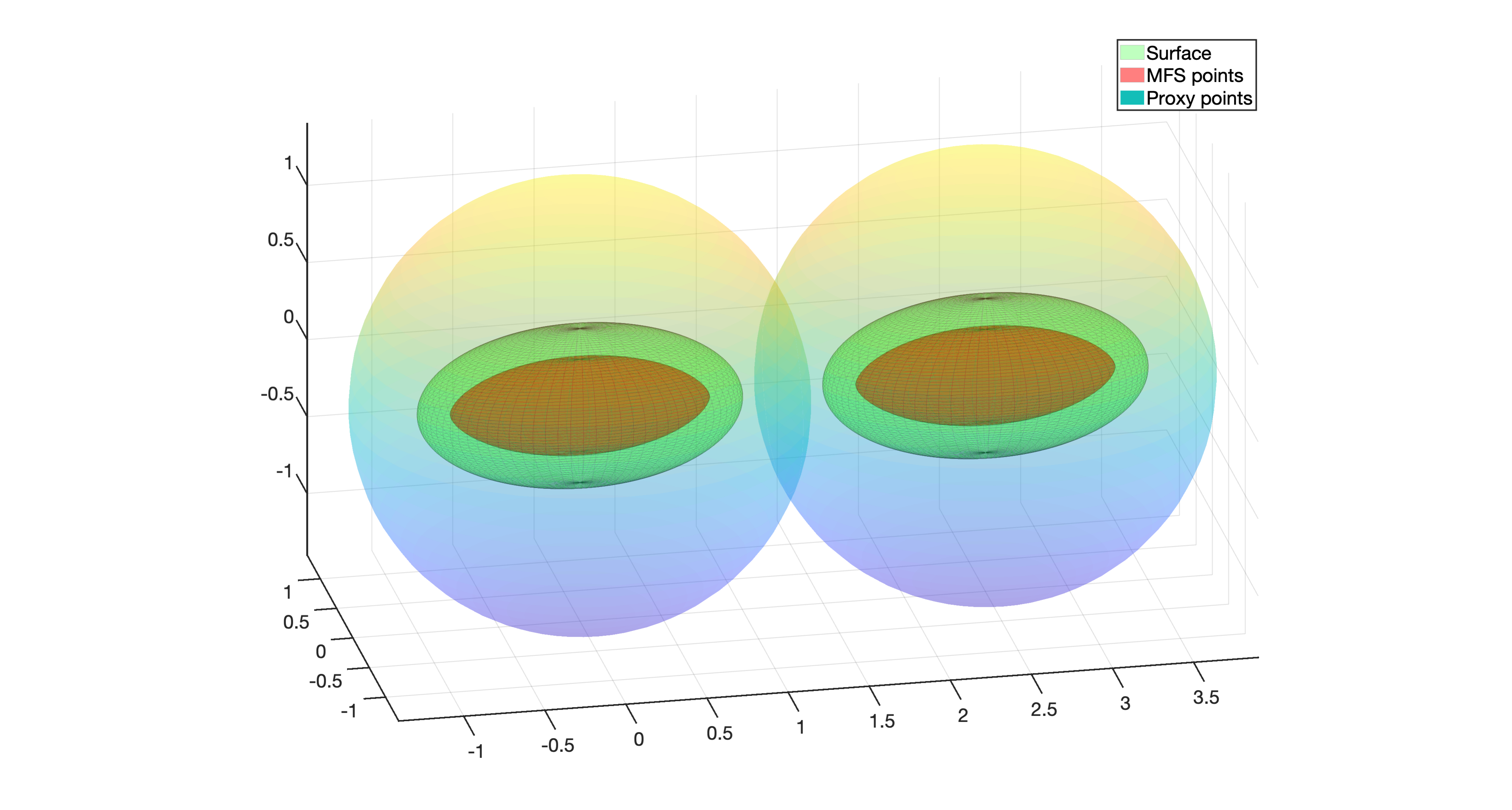}
        \caption{3D Ellipsoid geometry}
        \label{fig: Ellip3D}
    \end{figure} 

    \begin{table}[ht]
        \caption{Two-scatterer 3D experiments with $\varepsilon = 10^{-6}$}
        \label{tab:EllipsoidTorus3D_5_10_3}
        \setlength{\tabcolsep}{4pt}
        \begin{tabular}{c@{\quad}c}
        \begin{minipage}{0.5\textwidth}\centering\small
        \textit{(a)} Two 3D ellipsoids ($\kappa = 5$)\\[5pt]
        \begin{tabular}{cccccc}
        \toprule
        $N$ &
        $d$ &
        $N_{\rm skel}$ &
        $K$ &
        $E_{\rm inc}$\\
        \midrule
        256  & 0.10 &  188    & 4.06  & $4.67\times 10^{-2}$\\
        1024 & 0.10 &  311   & 3.01  & $2.43\times 10^{-3}$\\
        4096 & 0.10 &  340   & 2.94  & $1.51\times 10^{-5}$\\
        16384 & 0.10 &  339   & 2.87  & $5.62\times 10^{-9}$\\
        \bottomrule
        \end{tabular}
        \end{minipage}
        &
        \begin{minipage}{0.5\textwidth}\centering\small
        \textit{(b)} Two 3D tori ($\kappa = 5$)\\[5pt]
        \begin{tabular}{cccccc}
        \toprule
        $N$ &
        $d$ &
        $N_{\rm skel}$ &
        $K$ &
        $E_{\rm inc}$\\
        \midrule
        256  & 0.18 &  198    & 976  & $2.07\times 10^{-1}$\\
        1024 & 0.18 &  397   & 45.71  & $6.53\times 10^{-3}$\\
        4096 & 0.18 &  574   & 3.80  & $5.73\times 10^{-6}$\\
        16384 & 0.18 &  589   & 3.81  & $3.58\times 10^{-7}$\\
        \bottomrule
        \end{tabular}
        \end{minipage}
        \tabularnewline
        \begin{minipage}{0.5\textwidth}\centering\small
        \vspace*{10pt}
        \textit{(c)} Two 3D ellipsoids ($\kappa = 10$)\\[5pt]
        \begin{tabular}{cccccc}
        \toprule
        $N$ &
        $d$ &
        $N_{\rm skel}$ &
        $K$ &
        $E_{\rm inc}$\\
        \midrule
        256  & 0.50 &  209    & 19900  & $1.52\times 10^{-1}$\\
        1024 & 0.30 &  278   & 4.11  & $1.68\times 10^{-4}$\\
        4096 & 0.30 &  425   & 3.82  & $2.49\times 10^{-6}$\\
        16384 & 0.10 &  425   & 3.72  & $7.62\times 10^{-8}$\\
        \bottomrule
        \end{tabular}
        \end{minipage}
        &
        \begin{minipage}{0.5\textwidth}\centering\small
        \vspace*{10pt}
        \textit{(d)} Two 3D tori ($\kappa = 10$)\\[5pt]
        \begin{tabular}{cccccc}
        \toprule
        $N$ &
        $d$ &
        $N_{\rm skel}$ &
        $K$ &
        $E_{\rm inc}$\\
        \midrule
        256  & 0.18 &  221    & 6670  & $3.53\times 10^{0}$\\
        1024 & 0.18 &  451   & 876  & $2.12\times 10^{-1}$\\
        4096 & 0.18 &  656   & 4.07  & $2.36\times 10^{-4}$\\
        16384 & 0.18 &  678   & 4.09  & $6.03\times 10^{-6}$\\
        \bottomrule
        \end{tabular}
        \end{minipage}
        \end{tabular}
    \end{table}    
\end{example}

\begin{example}
\label{ex:scaling-3D}
    For our final experiment, we consider a scaling problem in 3D.
    The geometry is a mixture of tori and ellipsoids. As in the 2D scaling problem (Example~\ref{ex:2d-results-scaling}), the scattering matrix is identical for scatterers of the same shape (modulo phase factors upon rotation), which means that we can precompute the skeletonization and scattering matrices for a single ellipsoid and tori separately. Thus, for each error level, only two scattering matrices actually need to be formulated. \cref{tab:scaling_3d_k10-4} shows the numerical results for different numbers of scatterers $T$ under two prescribed error levels. Each column here shows the same quantity as in \cref{tab:scaling_2d_final}.

    The table shows that the number of matvecs grows only moderately (at most linearly) with the number of scatterers, and the time per matvec is also asymptotically linear (due to the use of the fast multipole method), which is comparable to the 2D scaling results. The local computations are also linear in the number of scatterers.
    
    The GMRES tolerance is here set equal to the target error level divided by 10. It can be seen that the solve time (and the number of matvecs) depends on the GMRES tolerance. If instead the same GMRES tolerance is used for both error levels, the number of matvecs will be almost independent of the error level (\textit{i.e.}, of the discretization resolution), as in the 2D scaling test (Example~\ref{ex:2d-results-scaling}, \cref{fig:gmres-tol-variation}).

\begin{table}[ht]
    \centering\small
    \caption{Scaling problem with 3D scatterers (mixture of ellipsoids and tori) and wavenumber $\kappa=10$} 
    \label{tab:scaling_3d_k10-4}
    \setlength{\tabcolsep}{2pt}
    \begin{tabular}{cccc|lll|l}
    \toprule
    \multicolumn{8}{c}{\textit{Error level $10^{-3}$ ($1024$ points per ellip., $4096$ points per tori, skel. tol.\ $10^{-3}$, }} \\
    \multicolumn{8}{c}{\textit{GMRES tol.\ $10^{-4}$, precomp. time \SI{16.67}{s})}} \\
    \midrule
    \multicolumn{4}{c|}{\bf Number of} & \multicolumn{3}{c|}{\bf Time [s]} & \multicolumn{1}{c}{\bf Error}\\
    scatterers &
    matvecs &
    points &
    skeletons &
    \multicolumn{1}{c}{local comp} &
    \multicolumn{1}{c}{solve} &
    \multicolumn{1}{c|}{per matvec} &
    \multicolumn{1}{c}{incoming} \\
    $T$ & $n_\text{matvec}$ & $N_\text{tot}$ & $N_\text{skel,tot}$ & \multicolumn{1}{c}{$t_\text{local}$} & \multicolumn{1}{c}{$t_\text{solve}$} & \multicolumn{1}{c|}{$t_\text{matvec}$} & \multicolumn{1}{c}{$E_\text{inc}$} \\
    \midrule
\num{4} & \num{20} & \num{10240} & \num{866} & \num[retain-zero-exponent=true]{2.35e+00} & \num{2.80e-01} & \num{1.40e-02} & \num{1.26e-03} \\
\num{8} & \num{26} & \num{20480} & \num{1732} & \num[retain-zero-exponent=true]{4.75e+00} & \num[retain-zero-exponent=true]{3.78e+00} & \num{1.45e-01} & \num{1.44e-03} \\
\num{16} & \num{30} & \num{40960} & \num{3464} & \num[retain-zero-exponent=true]{9.50e+00} & \num{1.76e+01} & \num{5.85e-01} & \num{1.41e-03} \\
\num{32} & \num{38} & \num{81920} & \num{6928} & \num{1.95e+01} & \num{2.62e+01} & \num{6.89e-01} & \num{1.68e-03} \\
\num{64} & \num{58} & \num{163840} & \num{13856} & \num{3.70e+01} & \num{4.69e+01} & \num{8.08e-01} & \num{2.63e-03} \\
\num{128} & \num{78} & \num{327680} & \num{27712} & \num{7.44e+01} & \num{2.19e+02} & \num[retain-zero-exponent=true]{2.81e+00} & \num{3.76e-03} \\
\num{256} & \num{118} & \num{655360} & \num{55424} & \num{1.47e+02} & \num{3.93e+02} & \num[retain-zero-exponent=true]{3.33e+00} & \num{3.83e-03} \\
\num{512} & \num{196} & \num{1310720} & \num{110848} & \num{2.95e+02} & \num{8.52e+02} & \num[retain-zero-exponent=true]{4.35e+00} & \num{4.79e-03} \\
    \bottomrule
    \toprule
    \multicolumn{8}{c}{\textit{Error level $10^{-6}$ ($5776$ points per ellip., $8464$ points per tori, skel. tol.\ $10^{-5}$, }} \\
    \multicolumn{8}{c}{\textit{GMRES tol.\ $10^{-7}$, precomp. time \SI{203.74}{s})}} \\    
    \midrule
    \multicolumn{4}{c|}{\bf Number of} & \multicolumn{3}{c|}{\bf Time [s]} & \multicolumn{1}{c}{\bf Error}\\
    scatterers &
    matvecs &
    points &
    skeletons &
    \multicolumn{1}{c}{local comp} &
    \multicolumn{1}{c}{solve} &
    \multicolumn{1}{c|}{per matvec} &
    \multicolumn{1}{c}{incoming} \\
    $T$ & $n_\text{matvec}$ & $N_\text{tot}$ & $N_\text{skel,tot}$ & \multicolumn{1}{c}{$t_\text{local}$} & \multicolumn{1}{c}{$t_\text{solve}$} & \multicolumn{1}{c|}{$t_\text{matvec}$} & \multicolumn{1}{c}{$E_\text{inc}$} \\
    \midrule
\num{4} & \num{30} & \num{28480} & \num{1644} & \num[retain-zero-exponent=true]{7.77e+00} & \num[retain-zero-exponent=true]{1.20e+00} & \num{4.00e-02} & \num{9.43e-07} \\
\num{8} & \num{38} & \num{56960} & \num{3288} & \num{1.62e+01} & \num{1.15e+01} & \num{3.02e-01} & \num{7.89e-07} \\
\num{16} & \num{46} & \num{113920} & \num{6576} & \num{3.24e+01} & \num{4.25e+01} & \num{9.24e-01} & \num{1.31e-06} \\
\num{32} & \num{64} & \num{227840} & \num{13152} & \num{6.77e+01} & \num{6.89e+01} & \num[retain-zero-exponent=true]{1.08e+00} & \num{1.46e-06} \\
\num{64} & \num{102} & \num{455680} & \num{26304} & \num{1.26e+02} & \num{1.22e+02} & \num[retain-zero-exponent=true]{1.19e+00} & \num{1.66e-06} \\
\num{128} & \num{138} & \num{911360} & \num{52608} & \num{2.54e+02} & \num{6.00e+02} & \num[retain-zero-exponent=true]{4.35e+00} & \num{2.47e-06} \\
\num{256} & \num{218} & \num{1822720} & \num{105216} & \num{5.00e+02} & \num{1.20e+03} & \num[retain-zero-exponent=true]{5.52e+00} & \num{3.01e-06} \\
\num{512} & \num{379} & \num{3645440} & \num{210432} & \num{1.04e+03} & \num{2.89e+03} & \num[retain-zero-exponent=true]{7.61e+00} & \num{2.98e-06} \\
    \bottomrule
    \end{tabular}
\end{table}

    As in the 2D scaling experiment, the observed error $E_\text{inc}$ grows slightly as the number of scatterers $T$ increases. The same remarks apply here, and the error can be kept fixed if the number of collocation points $N$ is increased slightly, and the GMRES tolerance reduced slightly. We emphasize that this is only a very marginal effect.
\end{example}

\section{Conclusions}

The paper describes a numerical method for solving acoustic multibody scattering problems in two and three dimensions. The solver is based on the method of fundamental solutions, and applies the ``skeletonization'' process of \cite{2005_martinsson_skel} to compress the number of unknowns, and to form numerical scattering matrices for each individual scatterer. A global system is then formed using the individual scattering matrices. A key finding is that this global system is relatively well-conditioned, despite the use of a numerically ill-conditioned local solver.

Numerical experiments demonstrate that the solver can even for complex geometries involving large number of scatterers compute highly accurate solutions. The solver was benchmarked against state-of-the-art BIE solvers, and attains similar accuracy, despite using a much simpler discretization technique. 

The method is designed to interact well with fast summation schemes such as the fast multiple method, as all interactions between scatterers occur through the free-space fundamental solution of the Helmholtz equation. 
We found that even for complex geometries involving trapping scatterers, the number of GMRES iterations tends to scale linearly in the number of scatterers. 

Importantly, since all scattering matrices are computed locally, the number of degrees of freedom in the global linear system is independent of the local geometric complexity. For instance, our numerical examples involve domains with sharp corners that require a large number of local degrees of freedom to fully resolve the local singularities. 

\section*{Declarations}

\bmhead{Conflict of interest}

The authors declare no competing interests.

\bmhead{Funding}

PGM is supported by
the Office of Naval Research (N00014-18-1-2354),
the National Science Foundation (DMS-2313434 and DMS-1952735), and
the Department of Energy ASCR (DE-SC0022251). JB is supported in part by the Peter O'Donnell Jr.\ Postdoctoral Fellowship at the Oden Institute.

\bmhead{Author contribution}

All authors read and approved the final manuscript. Yunhui Cai: validation, formal analysis, investigation, numerical experiments, visualization, writing -- original draft of \cref{sec:MFS,sec:constructS,sec:NumericalResults}, review and editing. Joar Bagge: validation, formal analysis, investigation, numerical experiments, writing -- original draft of \cref{sec:introduction,sec:NumericalResults}, review and editing. Per-Gunnar Martinsson: conceptualization, methodology, writing -- original draft of \cref{sec:introduction,sec:overview}, review and editing, supervision.

\bmhead{Acknowledgements}

The authors wish to thank Abinand Gopal of the University of California at Davis
for sharing important insights about the methods proposed.

\bmhead{Data availability}

No data was used for the research described in the article. 

\bibliography{references}

\end{document}